\documentclass{jfl}
\usepackage{array}
\usepackage{url}
\usepackage{amsmath}
\usepackage{amssymb}
\usepackage{amsfonts}
\usepackage{amsthm}
\usepackage{stmaryrd}
\usepackage{prooftree}

\def\causes{\triangleright}

\def\causRel#1{\mathrel{\mathbf{R}_{#1}}}
\def\entailsACR{\mathbin{\entails_{\text{ACR}}}}

\def\entails{\vdash}
\def\entailsRanked#1#2#3{\entails^{\!\!\raise 0.5ex 
    \hbox{$\scriptstyle#1,\,#2$}}_{\!\!#3}}
\def\rank{\textrm{rk}}

\def\models{\Vdash}
\def\forces{\Vdash}
\def\modelsNec{\models_{\nec}}
\def\eqdef{\overset{\text{\tiny def}}{=}}
\def\closure#1{\overline{#1}}

\def\lang{\mathfrak{L}}

\def\langMod{\lang_{\nec}}
\def\entailsMod{\entails_{\!\!\!\scriptscriptstyle\nec}}
\def\allModels{\mathfrak{M}}
\def\kripke#1{K_{#1}}
\def\semval#1{\llbracket#1\rrbracket}
\def\turnerMod{\textsf{\textbf{C}}}
\def\semEntMod{\vDash_{\nec}}

\def\name#1{\mathord{\scriptstyle\textsf{#1}}}
\def\shoot{\name{shoot}}
\def\wait{\name{wait}}
\def\alive{\name{alive}}
\def\loaded{\name{loaded}}

\def\nec{\oblong}

\def\bot{\perp}
\def\nat{\mathbb{N}}
\def\lequiv{\cong}
\def\ltrue{\top}
\def\lfalse{\perp}

\title{Causality, Modality, and Explanation}
\author[White]{Graham White}
\address{Department of Computer Science\\
Queen Mary, University of London\\
London E2 9QD}
\email{graham@dcs.qmul.ac.uk}
\urladdress{http://www.dcs.qmul.ac.uk/\urltilde graham}
\primarydata{03B45}
\secondarydata{03B45}
\secondarydata{03F05}
\begin{document}
\begin{abstract}
  We start with Fodor's critique of cognitive science in
  \cite{FodorJA:mindww}: he argues that much mental activity
  cannot be handled by the current methods of cognitive science
  because it is nonmonotonic and, therefore, is global in nature,
  is not context-free, and is thus not capable of being
  formalised by a Turing-like mental architecture. We look at
  the use of non-monotonic logic in the Artificial Intelligence
  community, particularly with the discussion of the so-called
  ``frame problem''. The mainstream approach to the frame problem
  is, we argue, probably susceptible to Fodor's critique: however,
  there is an alternative approach, due to McCain and Turner, which
  is, when suitably reformulated, not susceptible. 
  In the course of our argument, we give a proof theory for
  the McCain-Turner system, and show that it satisfies cut
  elimination. We have two substantive conclusions: firstly, that
  Fodor's argument depends on assumptions about logical form which
  not all non-monotonic theories satisfy; and, secondly, that 
  metatheory plays an important role in the context of evolutionary
  accounts of rationality. 
\end{abstract}
\maketitle
\tableofcontents
\setcounter{section}{-1}
\section{Introduction}
Fodor has argued \cite{FodorJA:mindww} that mental processes fall into
two classes: those which are modular (that is, effectively
encapsulated from other mental processes), and those which are
not. The modular ones take place by means of local, syntactic
operations on mental representations \cite[pp.~18f.]{FodorJA:mindww}:
if they function properly, they will be insensitive to what may or may
not happen outside their module.  In fact, more than this is true:
such processes are \emph{context invariant}
\cite[pp.~34f.]{FodorJA:mindww}.  All they care about is the syntax of
the mental representations which they are working on. In particular,
they are monotonic: the validity of the inferences they perform is
unaffected by adding extra mental representations to their context.

The non-modular ones, by contrast, are not monotonic. Fodor's examples
of such non-modular processes are, firstly, abduction, and, secondly,
the belief change apparatus associated with abductive inference
(including judgements as to the centrality of particular beliefs in
one's belief system). These processes are, plausibly, non-monotonic
and context-sensitive \cite[pp.~33ff.]{FodorJA:mindww}: and thus,
Fodor goes on to argue, they cannot be described as operations on the
(context-invariant) syntax of mental representations.

Now one \emph{could} quibble about all this in various ways: one could
worry about whether Fodor had distinguished sufficiently well between
abstract and concrete syntax
\cite{mccarthy62:_towar_abstr_scien_of_comput}, one could wonder
whether he had taken sufficiently into account the possibility of
context dependence which was regular enough not to have the
deleterious consequences that he wishes to avoid
\cite{strachey67:_fundam_concep_of_progr_languag} -- but, for all
that, his argument is cogent and worrying.

Rather than raising these sort of objections to Fodor's argument, it
may be more illuminating to look at how the artificial intelligence
community deals with this sort of problem. We will look at what is
described in the AI community as the ``frame problem'': although this
problem, which is mainly concerned with reasoning about change in the
world, is somewhat different from the belief change problems that
Fodor discusses, the two areas are systematically related
\cite{bochman01:_logic_theor_of_nonmon_infer}, so that discussions of
one should, \emph{mutatis mutandis}, apply to the other.

Our results are interesting. We investigate two approaches to the
frame problem: one of them, the more mainstream one, is based on the
idea of minimising change, and probably \emph{is} susceptible to
critique in Fodor's manner. However, the other one, due to McCain and
Turner, is probably not thus susceptible: it shows illuminatingly how
reasoning can be nonmonotonic without being catastrophically
global.

In arguing for this position, we have work, both technical and
philosophical, to do.  The McCain-Turner system is usually presented
model-theoretically: we develop a proof theory for it, in order to
show directly what reasoning in that system is like.  In particular,
we prove a cut elimination result, which can be regarded as showing
how the system can support nonmonotonic inference and still be
tractable. Proofs of these results are quite technical, and have been
modularised into the appendix of this paper.

We also have philosophical work to do. McCain and Turner present their
system as a system of what they term ``causal reasoning'': I argue
that the system can, much more plausibly, be presented as a logic of
\emph{explanation}. This, of course, brings it much closer to Fodor's
concern with large-scale cognitive architecture, in which explanation
ought to play a central role. So, in order to make this claim more
plausible, I describe an interesting logic of argument due to Parsons
and Jennings: this logic turns out to be merely a notational variant
of a fragment of the McCain-Turner system.

\section{The Logical Background}
\subsection{Basics}
Regardless of the details of the common approaches to reasoning about
action, they all have several things in common. They all assume that
one should start with a logical description of actions and their
effects: they all (or mostly all) assume that the logical apparatus
should deliver a description of the sequence of states of affairs
resulting from a given sequence of actions.  We will, for the sake of
simplicity, assume that we have linear discrete time, indexed by the
natural numbers: propositions will have an argument place for the
time, into which terms standing for temporal instants can be
substituted. Devotees of more sophisticated approaches -- for example,
of Reiter's elegant situation calculus formulation
\cite{reiter01:_knowl_action} -- are welcome to translate my
oversimplified notation into their variant.

\subsection{The Frame Problem}
A fundamental difficulty behind using logic to reason about action in
this way is what is known as the frame problem. The problem is that,
given the state of the world before an action, there are many, merely
``logically possible'', states of the world after it; but most of
these states are physically unrealistic. This remains so even if we
give what seem to be perfectly adequate logical formalisations of
actions.

Here is a famous example of what is at issue. 
\begin{example}[The Yale Shooting Problem]\label{example:ysp}
  Suppose that we describe a state of the world at a particular time
  $t$ by propositions $\alive$ and $\loaded$. The world (a
  particularly simple one) contains a victim, a gun (and, presumably,
  a shooter, but he or she can remain unformalised).  $\alive_{t}$ will
  say that some victim is still alive at $t$, whereas $\loaded_{t}$
  will say that a gun is loaded at that time. Consider actions
  $\shoot$ and $\wait$: $\wait$ has no effect, whereas $\shoot$, if
  the gun is loaded, kills the victim. So we can specify the shoot action
  as follows (for the sake of simplicity, we assume that every action
  takes one unit of time):
  \begin{equation}
    \shoot_{t} \land \loaded_{t} \entails \lnot\alive_{t+1}
  \end{equation}
  Suppose now we have the following scenario: the
  gun is loaded at time $0$, then there is a $\wait$ action,
  and the gun is then shot. So the following sequence of truth values
  is compatible with our specification:
  \begin{equation}
    \begin{array}{c|ccc}
     t & 0 & 1 & 2 \\
       &\wait & \shoot & \\
       \hline
       \alive & \ltrue & \ltrue & \lfalse\\
       \loaded & \ltrue &\ltrue & \lfalse\\
       \hline
    \end{array}
    \label{eq:goodHistory}
  \end{equation}
  Unfortunately, so is this:
  \begin{equation}
    \begin{array}{c|ccc}
     t & 0 & 1 & 2 \\
       &\wait & \shoot & \\
       \hline
       \alive & \ltrue & \ltrue & \ltrue\\
       \loaded & \ltrue &\lfalse & \lfalse\\
       \hline
    \end{array}\label{eq:badHistory}
  \end{equation}
\end{example}

However one might wish to diagnose the situation here, it
is clear that something has gone wrong: logic, supplemented
by unproblematic-seeming axiomatisations of the effects of actions,
generates more solutions to these problems than it should. 
We need extra constraints, and there have been at least
two attempts to say what sort of constraints we should have:
a mainstream approach, based on the idea of minimisation, 
and an alternative approach, rather less clearly formulated, but
technically more viable and more interesting. As I shall argue,
the alternative approach is, in fact, based on the idea of
explanation.  

\subsection{Minimisation-Based Approaches} 
We first outline, and criticise, the minimisation-based approaches: 
this will not only show some of the issues at stake, but it will
also give an example of an approach which is, arguably, susceptible
to Fodor's critique. 

The basic idea behind these approaches is simple: that we should pick
those histories which minimise the amount of change. One can,
presumably, motivate this criterion as follows: the axiomatisations of
actions give the changes which \emph{must} take place if the actions
are to be performed, and one wants exactly these changes and no
others.  Hence one should have a minimal amount of change. This idea
is very close to Lewis' treatment of counterfactuals
\cite{00LewisD:cau,00LewisD:coucp,00LewisD:coudta}, and thus has a
respectable philosophical pedigree. It can, furthermore, be supported
by a mathematical argument: a wide class of non-monotonic logics can
be specified by giving a preference relation on models of some
monotonic theory, and then selecting only those models which are
minimal according to the preference relation \cite{brewka}. So, this
approach seems to have the advantage of a strong philosophical
motivation, together with mathematical generality.

However, problems arise when one tries to make this approach more
specific. Example~\ref{example:ysp} shows the first of these problems:
there are two histories in that example, and the wrong one has one
change of truth value, while the correct one has two changes. Even if
we minimise by looking at containment relations between sets of
changes, rather than naively counting them, the problem persists. So
minimisation often gives incorrect results. One can, perhaps, fix it,
but many of the fixes (for example, preferring later changes to
earlier ones) seem to be merely \emph{ad hoc}.

The second problem is this: what does one mean by ``amount of
change''?  This is, admittedly, a vague description, but it is not
clear how to make it any less vague. The usual approach is this: one
fixes a language, and then, for each history, one fixes the set of
truth-value changes (that is, the set of ordered pairs consisting of a
time and a primitive of the language which changes truth-value at that
time). Containment between these sets gives a preference relation
between histories: those histories are to be preferred which have
minimal sets of truth-value changes. One can give this procedure a
much more sophisticated and formal definition: it is generally known
as \emph{circumscription} \cite{00LifschitzV:cir}, and is widely
used in the AI community.

The problem with this approach is that it is strongly
language-dependent: it is possible to have equivalent languages, with
different primitives, which give different preference relations
according to this definition. One can, by a suitable choice of
primitives for the language, make the preference relation come out any
way that one wants \cite{WhiteGG:intc}. In more technical language,
the entailment relations which circumscription yields are not
preserved by uniform substitution.  Failure of uniform substitution is
(\emph{pace} \citend{makinson03:_ways_doing_logic}) a cause of concern,
for two reasons.  Practically, it makes the use of these theories
rather problematic: there may well be some magical set of primitives
which make the predictions come out right, but one is given no
guidance as to what it might be. And, theoretically, this remains a
worry: closure under uniform substitution is motivated by a concern
that mathematics should be more than mere syntax, that is, that the
entities which mathematics and logic describe should be independent of
our choice of the language -- the ``coordinate frame'', as it were --
which we use to describe them.

Now Lewis works with possible worlds, which, for him, are first class
individuals, and his closeness relation between possible worlds --
which corresponds to our notion of minimising change -- is simply
given: he does not say how to work it out, and, in particular, how to
work it out on the basis of physically realistic
measurements. Consequently, he does not explicitly face the problems
that we have described above. However, we should note that his account
may well be susceptible to similar problems: it may work in the
abstract, but, if we require that the nearness relations between
worlds should be capable of evaluation in terms of continuously
differentiable functions of measurements made in those worlds, we may
end up with a theory incompatible with physics
\cite{white00:_lewis_causal_possib_world}. Lewis' approach turns out
to be, in Fodor's terms, \emph{global}: we have to know about all of
the possible worlds that there are, and we have to know these worlds
directly as individuals, in order to make causal inferences.

The difficulty with Lewis' account comes from a quite general result
(specifically, from what is known as Noether's theorem) which has to
do with the invariance of the laws of physics under coordinate
transformation. So, conceptually, these considerations are not
infinitely far removed from the problems that circumscription has with
invariance under uniform substitution: and the problems with
this approach seem to be quite deeply rooted. 

It is, of course, not impossible that one could modify circumscription so
as to allay these worries, or, indeed, that Lewis' account of 
counterfactuals could be made specific enough to work with physical
laws. However, these are difficult problems, and they may well
be better illuminated by an alternative approach. 

\subsection{Explanation-Based Approaches} 

Consider the bad history given in \eqref{eq:badHistory}. We
might attempt a critique of it as follows: the problem is not
that there are \emph{too many} changes, but that we cannot
\emph{explain} one of the facts -- namely, the falsity of
$\loaded$ at time 1. All of the other facts in that
history can be explained, either because of the effects of actions
or because they are facts unchanged from the previous time;
but not this one. So a good constraint to put on our logic
might be to require that every fact should be explained. 

This approach might have a good deal to recommend it. A good notion of
explanation would, one thinks, be language-, or coordinate-frame-,
independent: the translation of an explanation would surely still be
an explanation. Furthermore, it has a pleasing directness: human
beings like us perform a lot of common sense causal explanation, and
reflection on this process might give us useful data for formalising
such explanations. We might also cite the philosophical tradition:
since the time of Aristotle, at least, causality and explanation have
been very closely linked \cite[p.~32]{barnes82:_arist},
\cite{sep-aristotle-causality}.

In order to carry this programme forward, we will have to have a
certain amount of technical apparatus. We will have to have a notion
of explanation: the things to be explained will be facts, and we
should be able to express, of any particular fact, whether it is
explained or not. These notions should be closed under uniform
substitution: that is, the translation of an explanation should be an
explanation.  Such an approach to reasoning about action is, indeed,
possible: it was formulated by McCain and Turner, and we will study it
in the remainder of this paper.

\section{McCain and Turners' System}
Their treatment \cite{McCainTurner:cautac} of the frame problem goes
as follows. Suppose that we have a language $\lang$, which will have
the connectives of the classical propositional
calculus.  McCain and Turner
then consider a binary sentential operator: let us write it $\cdot
\causes \cdot$.  Applications of this operator -- of the form $\phi
\causes \psi$ -- are called ``causal rules'', and they can, in McCain
and Turner's treatment, be regarded as purely metatheoretic
assertions: \emph{they} read them as `$\phi$ causes $\psi$', but we will
argue for a reading of `$\phi$ explains $\psi$'.  A collection of
causal rules will be called a \emph{causal theory}.  Given a causal
theory $\Theta$ in our language $\lang$, we define an operator
$(\cdot)^{\Theta}$, from models of $\lang$ to sets of sentences, as
follows:
\begin{equation}
  \label{eq:causExpModel}
  M^{\Theta} \quad = \quad \closure{\{ \psi | 
\phi \causes \psi \in D, M \models \phi\}},
\end{equation}
where $\closure{\cdot}$ is closure under entailment. 
We now say that a model $M$ of $\lang$ is \emph{causally explained}
iff it is the only model of $M^{\Theta}$.  And we say that a proposition
$P$ is a \emph{consequence} of $\Theta$ iff it is true in all of the
$\Theta$-causally explained models of $\lang$. 
\subsection{McCain and Turner's Causal Rules}%
\label{sec:mccainTurnerRules}

We will give a simplified version of McCain and Turner's causal
theory: in order to make our exposition more direct, we will be
presenting a version of our sequent calculus without the quantifiers
(we could describe a version with quantifiers, but it would add extra
technical complexity, and this extra complexity would not be germane
to the main argument of this paper).

Theories of this sort will talk of two sorts of entities,
\emph{fluents} and \emph{actions}: a fluent is something which can be
true or false at a particular time, whereas an action is something
which can be performed (or not) at a particular time. We represent
each of these by a series of temporally indexed propositional atoms:
\begin{align}
  f_t&\quad (t \in \nat) &\quad\text{for each fluent} \\
  a_t&\quad(t \in \nat)  &\quad\text{for each action}
\end{align}
Intuitively, $f_t$ is true iff $f$ holds at time $t$, and $a_t$ is
true iff the action $a$ occurs at time $t$.  The atoms $f_t$ will be
called the \emph{fluent atoms}: a \emph{fluent literal} will be a
fluent atom or the negation of a fluent atom.

We now describe our causal theory. 
We assume that the effects of the actions can be formalised by giving,
for each action $a$, a finite set of \emph{precondition-postcondition
  pairs}. Such a pair consists of a precondition and a postcondition,
each of them a conjunction of fluent literals: the intuitive meaning 
of it is that, if the precondition is true in a situation,
then, after the execution of the relevant action, the postcondition
will be true in the successor situation. For example, the
precondition of the $\shoot$ action is the $\loaded$ fluent,
and the postcondition of the $\shoot$ action is $\lnot \alive
\land \lnot \loaded$. 
 
So, for each action $a$, and for each precondition-postcondition pair
$\langle f,g \rangle$ belonging to that action, we add a series of causal
rules: 
\begin{align}
  (f_t \land a_t) &\causes g_{t+1} 
  \intertext{where $f$ is the precondition of $a$, and $g$ is
    the postcondition. We also have to say which
    actions occur at which times: so, if action $a$ occurs at time
    $t_0$, we add the causal rule} 
  a_{t_0}&\causes a_{t_0}
  \intertext{We also need to deal with fluents that remain unchanged
    (see \cite{ShanahanM:solFP}): so, for every fluent literal $f_t$,
    we add the causal rule} 
  f_t \land f_{t+1} &\causes  f_{t+1}\label{eq:persistenceRule} 
  \intertext{Finally, we assume that
    a set of fluents is given which describes the initial situation:
    for each such fluent $f$, we add the causal rule} 
  f_0  &\causes f_0.\label{eq:initialStateRule}
\end{align}
Given these rules, this procedure appears to work: that is, it yields
correct solutions of the frame problem (the reader is invited to check
that it works with Example~\ref{example:ysp}). As a bonus, it is
fairly efficient computationally.
\subsection{The Road Ahead}
As it stands, though, McCain and Turner's definitions are not entirely
perspicuous, either mathematically or conceptually. The definition of
causally explained models is mysterious: we would, thus, like a more
illuminating treatment of the mathematics. Furthermore, it is
difficult to explain the conceptual role of $\causes$. For example,
one of McCain and Turner's rules is $ f_0
\causes f_0$ \eqref{eq:initialStateRule}. 
But reading $\causes$ as `causes' here is simply
implausible: propositions -- or the states of affairs which they
denote -- are not usually thought of as causing themselves.

So, we will do two things in this paper. Firstly, we will give an
alternative definition of \eqref{eq:causExpModel}: it can easily be
reformulated as a relation between models, and relations on a set of
models give modal operators. So we can, instead, define a suitable
modal operator, $\nec$; given this operator, we can reformulate 
the definition of $M^{\Theta}$ 
\eqref{eq:causExpModel} as 
\begin{align*}
  M^{\Theta} \quad &= \quad \{ \phi | M \models \nec \phi \};
  \intertext{similarly, the relation}
  \phi &\entails \nec \psi
  \intertext{is equivalent to}
  \phi &\causes \psi
\end{align*}
(equivalent in the sense that the
two relations give the same set of ``causally closed''
models).

But, as well as merely technical reformulations, we also want to say
what these constructions \emph{mean}.  We argue, then, that a better
reading of $\phi \causes \psi$ is that `$\phi$ explains $\psi$' --
and, correspondingly, $\nec \psi$ can be regarded as a disjunction of
all of the possible explanations of $\psi$. In fact, this explanatory
reading has been anticipated: \citend[p.~451]{LifschitzV:logce}
paraphrases $\phi \causes \psi$ as ``[$\psi$] has a cause if [$\phi$]
is true, or \ldots\ [$\phi$] provides a `causal explanation' for
[$\psi$]''. And, given this explanatory reading, reflexive rules seem
far less contentious: \eqref{eq:causExpModel} simply says that the
proposition $f_0$ ``explains itself''. Self-explanation is a good
deal less problematic than self-causation, since explanations, after
all, have to come to an end (or, of course, loop) at some point --
see, for example \citend[\S217]{wittgenstein67:_philos_inves}. In a
similar way, McCain and Turner's rule for persistence
\eqref{eq:persistenceRule} can be easily motivated in terms of
explanation: a good explanation for a fluent being true at $t+1$ is
that it was true at $t$, and that its truth value is unchanged between
$t$ and $t+1$. But it is hard to read this rule \emph{causally}
without, again, invoking self-causation.

\section{The Modal Logic}
\subsection{Motivation}
Our system is defined proof-theoretically, but it was arrived at by
reformulating McCain and Turner's original definition: our aim was to
find a sequent calculus, with a good proof theory, which corresponded
very closely to their model-theoretic construction. We proceed by
progressively reformulating their work.

\begin{remark}[Notation]
  We should explain a few notational conventions. $\phi$ and $\psi$ will always
  stand for, respectively, the antecedents and consequents of causal rules. General
  propositions will be written using lower case Roman letters, $p, q, r, \ldots$. 
  Sets of propositions will be written with upper case Greek letters, $\Gamma, \Delta$;
  $\Theta$ will be used for causal theories, i.e.\ sets of causal rules. 
  Languages will be written $\lang$, with subscripts and superscripts.
  Models of our various theories will be written $M, M', \ldots$ And
  when, in the appendix, we do ordinal analysis we will use lower case Greek letters,
  but not $\phi$ or $\psi$.
\end{remark}

\subsubsection{Reformulation 1: Model Theory}
First a definition:
\begin{definition}\label{def:causRel}
  A causal theory $\Theta$ defines a relation $\causRel{\Theta}$ on
  models of the language $\lang$ by
\begin{equation}
  \begin{split}
    M \causRel{\Theta} M' \: \Leftrightarrow \:&
    \text{for every } \phi \causes \psi \in \Theta,\\
    &\text{if } M \models \phi, \text{ then } M' \models \psi
  \end{split}
\end{equation}
\end{definition}
We then have:
\begin{lemma}
  For any $p \in \lang$, 
  \begin{equation}
    M^{\Theta} \entails p \:\Leftrightarrow \:
    M' \models p \text{ for all } M' \text{ with } M \causRel{\Theta} M'
  \end{equation}
\end{lemma}
\begin{proof}
  Note first that, for any model $M'$,  
  \begin{displaymath}
    M' \models M^{\Theta} \: \Leftrightarrow\:
    \text{for all } \phi \causes \psi \in \Theta. \, 
    M \models \phi \Rightarrow M' \models \psi,
  \end{displaymath}
  so we have 
  \begin{displaymath}
    \begin{split}
      M^{\Theta} \entails p \:&\Leftrightarrow \:
      M' \models p \text{ for all } M' \models M^{\Theta}\\
      &\Leftrightarrow \:\text{for all } M' 
      \begin{aligned}[t]
        (\text{for all}&\: \phi \causes \psi \in \Theta\; 
        M \models \phi \Rightarrow M' \models \psi) \\
        &\Rightarrow M' \models p
      \end{aligned}\\
      &\Leftrightarrow \:\text{for all } M'. 
      M \causRel{\Theta} M' \Rightarrow M' \models p
    \end{split}
  \end{displaymath}
\end{proof}
So we have 
\begin{corollary}
  For any model $M$, $M$ is causally explained according to
  $\Theta$ iff
  \begin{displaymath}
    \{M\} \: = \: \{ M' | M \causRel{\Theta} M'\}
  \end{displaymath}
\end{corollary}

\subsubsection{Reformulation 2: Modal Models}\label{kripkeSection}
We can now reformulate McCain and Turner's theory in modal terms.
\begin{definition}\label{def:canMod}
  Given, as above, a language $\lang$ and a causal theory
  $\Theta$, define a Kripke model $\kripke{\Theta}$ as follows:
  \begin{description}
  \item[the language] is $\langMod$, generated by $\lang$ together with
    the modal operator $\nec$.
  \item[the frame]consists of the set $\allModels$ of models of the
    non-modal language $\lang$, together with the accessibility
    relation $\causRel{\Theta}$.
  \item[the forcing relation]is given by the usual $\models$ relation
    between elements of $\allModels$ and propositions of $\lang$,
    extended to modal formulae in the usual way. If we wish to be
    pedantic (but we rarely will), we could call the new forcing
    relation $\models_{\nec}$.
  \end{description}
\end{definition}
Notice that we have:
\begin{proposition}\label{modalCausalAxiomProp}
  For each causal law $P \causes Q$, the modal sentence 
  \begin{equation}
    P \rightarrow \nec Q\label{causalModProp}
  \end{equation}
  is valid at each world of $K$. 
\end{proposition}
\begin{proof}
  This follows directly from the definition of $\causRel{\Theta}$.
\end{proof}
Then we have:
\begin{proposition}
  For $M \in \allModels$, $M$ is causally explained by $\Theta$ iff,
  as a world of the Kripke model $\kripke{\Theta}$,
  \begin{equation}
    M \models (\nec p \rightarrow p) \land 
    (p \rightarrow \nec p)\label{modalCrit}
  \end{equation}
for any $p$. 
\end{proposition}
\begin{proof}
  This is a standard modal reformulation of the condition
  \begin{displaymath}
    \{M\} = \{M' | M \causRel{\Theta} M'\}; 
  \end{displaymath}
  see, for example, \citend{vanBenthemJ:hanplii}.
\end{proof}
This result is pretty rather than useful: it refers to a particular
Kripke model, and models are hard to present in any effective sense.

We can make progress, though, by asking how \emph{other} Kripke models
are related to $\kripke{\Theta}$. 
\begin{definition}\label{definition:finalInterpretation}
  Let $K$ be any other Kripke model of $\langMod$
  such that, for any $\phi \causes
  \psi \in \Theta$, and any world $w$ of $K$,
  \begin{displaymath}
    w \forces \phi \rightarrow \nec \psi. 
  \end{displaymath}
  Then define a map 
  \begin{displaymath}
    \eta: \text{worlds of}\: K \:\: \rightarrow \:\: 
    \text{worlds of}\:\kripke{\Theta}
  \end{displaymath}
  by sending a world $w$ of $K$ to
  \begin{equation}
    \{ p\, |\, p \in \lang, w \forces p\} \label{eq:4}
  \end{equation}
  Note that the sets of propositions defined by (\ref{eq:4}) are, by
  the properties of Kripke models, models of $\lang$: consequently,
  these sets of propositions are, in fact, worlds of $\kripke{\Theta}$.
\end{definition}

\begin{proposition}
  For $w, w'$ worlds of $K$, if $w R_{K} w'$,
  then $\eta(w) R_{\Theta} \eta(w')$. 
\end{proposition}
\begin{proof}
  Suppose that $\phi \causes \psi \in \Theta$,
  and that $\eta(w) \forces \phi$. Then, by definition
  of $\eta$, $w \forces \phi$. However, by the assumption
  on $K$, $w \forces \phi \rightarrow \nec \psi$: since
  $w R_{K} w'$, we have $w' \forces \psi$, and so,
  by the definition of $\eta$, $\eta(w') \forces \psi$.
  This is so for any $\phi \causes \psi \in \Theta$,
  and thus we have $\eta(w) R_{\Theta} \eta(w')$.
\end{proof}
\begin{proposition}\label{proposition:modelMinMod}
  For any world $w$ of $K$, and any non-modal proposition
  $p$, if $\eta(w) \forces \nec p$, then $w \forces \nec p$.
\end{proposition}
\begin{proof}
  Immediate from the above. 
\end{proof}

This gives us an intuition as to what characterises the model
$\kripke{\Theta}$: among all of the models which globally satisfy
(\ref{causalModProp}), it is the one in which $\nec$ is strongest (or,
alternatively, the one in which the relation $R$ is the largest). This
shows us how to find a sequent calculus which corresponds to the model
$\kripke{\Theta}$, and which, consequently, represents the McCain-Turner
procedure.

There is a standard way of getting such a sequent calculus (see, for
example, \citend{hallnaes90:_proof_theor_approac_logic_progr,%
  hallnaes91:_proof_theor_approac_logic_progr}). We first give
ourselves right rules for $\nec$ corresponding to the constraints
which $\nec$ must satisfy: then we find left rules which invert the
right rules. The constraints which $\nec$ must satisfy are quite
simple: $\phi \rightarrow \nec \psi$ must be a theorem, for every
$\phi \causes \psi \in \Theta$. Rather more generally, if we have causal
rules $\{ \phi_i \causes \psi_i\}_{i \in I}$, and if we have
$\{\psi_i\}_{i \in I} \entails p$, then we must have $\{\phi_i \}_{i
  \in I} \entails \nec p$: this gives us a right rule for $\nec $
which ensures that $\nec$ is a \textbf{K} modality. The left rule will,
correspondingly, invert this right rule: and we get the rules for
$\nec$ in Table~\ref{rulesTable}. And, as the remainder of this paper
will show, our system does, indeed, give a sequent calculus for the
McCain-Turner procedure. It allows us to recover the Kripke model
$\kripke{\Theta}$, which turns out to be the canonical model of our modal
logic.

\subsection{The System}
\noindent Our system will be given by a sequent calculus,
as in \citend{WhiteGG:modfmt}; it is given by the rules in
Table~\ref{rulesTable}, whose formulation depends on an underlying set
of causal rules. As explained above, we will, for the
sake of simplicity, present the propositional version
of this system: we could present a first-order system,
but at the price of extra technical complexity. 

We should note that this sequent calculus introduces a new 
language, $\langMod$ -- an extension of $\lang$ by the
modal operator $\nec$ -- and a new consequence relation,
$\entailsMod$. Our cut elimination theorem will show that 
$\langMod$ is a definitional extension of $\lang$, and that
$\entails$ is the restriction of $\entailsMod$ to $\lang$: but,
for the moment, we will be careful to respect the differences
between the two entailment relations. 
\begin{table}
  \begin{center}
    \rule{\linewidth}{\proofrulebreadth}\\
    \caption{Sequent Calculus Rules}\label{rulesTable}
    \begin{displaymath}
      \setlength{\extrarowheight}{15pt}
      \addtolength{\arraycolsep}{12pt}
      \begin{array}{ll}
        \multicolumn{2}{c}{
          \begin{prooftree}
            \justifies
            p \entailsMod p
            \using
            \text{Ax}
          \end{prooftree}
        }
        \\
        \begin{prooftree}
          \justifies
          \lfalse \entailsMod
          \using
          \lfalse\text{L}
        \end{prooftree}
        &
        \begin{prooftree}
          \justifies
          \entailsMod \ltrue
          \using
          \ltrue\text{R}
        \end{prooftree}
        \\
        \begin{prooftree}
          \Gamma \entailsMod \Delta
          \justifies
          \Gamma, p \entailsMod \Delta
          \using
          \text{LW}
        \end{prooftree}&
        \begin{prooftree}
          \Gamma \entailsMod \Delta
          \justifies
          \Gamma \entailsMod p, \Delta
          \using
          \text{RW}
        \end{prooftree}
        \\
        \begin{prooftree}
          \Gamma, p, p \entailsMod \Delta
          \justifies
          \Gamma, p \entailsMod \Delta
          \using
          \text{LC}
        \end{prooftree}
        &
        \begin{prooftree}
          \Gamma \entailsMod p, p, \Delta
          \justifies
          \Gamma \entailsMod p, \Delta
          \using
          \text{RC}
        \end{prooftree}
        \\
        \begin{prooftree}
          \Gamma \entailsMod p, \Delta
          \justifies
          \Gamma, \lnot p \entailsMod \Delta
          \using
          \lnot \text{L}
        \end{prooftree}
        &
        \begin{prooftree}
          \Gamma, p \entailsMod \Delta
          \justifies
          \Gamma \entailsMod \lnot p, \Delta
          \using
          \lnot \text{R}
        \end{prooftree}
        \\
        \begin{prooftree}
          \Gamma, p, q \entailsMod \Delta
          \justifies
          \Gamma, p \land q \entailsMod \Delta
          \using
          \land \text{L}
        \end{prooftree}
        &
        \begin{prooftree}
          \Gamma \entailsMod p, \Delta \quad \Gamma \entailsMod q, \Delta
          \justifies
          \Gamma \entailsMod p \land q, \Delta
          \using
          \land \text{R}
        \end{prooftree}
        \\
        \begin{prooftree}
          \Gamma, p \entailsMod \Delta
          \quad
          \Gamma, q \entailsMod \Delta
          \justifies
          \Gamma, p \lor p \entailsMod \Delta
          \using
          \lor \text{L}
        \end{prooftree}
        &
        \begin{prooftree}
          \Gamma \entailsMod p, q, \Delta
          \justifies
          \Gamma \entailsMod p \lor q, \Delta
          \using
          \lor \text{R}
        \end{prooftree}
        \\
        \begin{prooftree}
          \Gamma \entailsMod p, \Delta 
          \quad
          \Gamma, q \entailsMod \Delta
          \justifies
          \Gamma, p \rightarrow q \entailsMod \Delta
          \using
          \rightarrow \text{L}
        \end{prooftree}
        &
        \begin{prooftree}
          \Gamma, p \entailsMod q, \Delta
          \justifies
          \Gamma \entailsMod p \rightarrow q, \Delta
          \using
          \rightarrow \text{R}
        \end{prooftree}
        \\
        \multicolumn{2}{c}{
          \begin{prooftree}
            \Gamma \entailsMod \phi_1\land \ldots \land\phi_k, \Delta
            \quad 
            \psi_1, \ldots, \psi_k \entailsMod p
            \justifies
            \Gamma \entailsMod \nec p, \Delta
            \using
            \nec \text{R}
          \end{prooftree}}
        \\
        \multicolumn{2}{c}{
          \begin{prooftree}
            \raise1pt\hbox{$
              \left\{
                \Gamma, \{\phi_i\}_{i \in I_j} \entailsMod
                \Delta,
                \quad 
                \{\psi_i\}_{i \in I_j}\entailsMod
                p
              \right\}_{j \in J}
              $}
            \justifies 
            \Gamma, \nec p \entailsMod
            \Delta \using \nec \text{L}
          \end{prooftree}}
        \\
        \multicolumn{2}{c}{
          \begin{prooftree}
            \Gamma \entailsMod p^m, \Delta \quad \Gamma', p^n 
            \entailsMod \Delta'
            \justifies
            \Gamma, \Gamma' \entailsMod \Delta, \Delta'
            \using
            \text{multicut}
          \end{prooftree}}
      \end{array}
    \end{displaymath}
    \begin{tabular}{rp{22pc}}
      \multicolumn{2}{l}{\textbf{conditions on the rules:}}\\
      $\nec \text{R}$
      &
      where,  for all 
      $i$, $\phi_i \causes \psi_i$ is
      a causal rule. 
      \\
      $\nec \text{L}$ 
      &
      for each
      appropriate $i$, we have 
      a causal rule $\phi_{i} \causes
      \psi_{i}$, and where the 
      $\{\phi_i\causes \psi_i\}_{i \in I_j}$ 
      are the only such finite sets of
      $\phi$s and $\psi$s, such that 
      $\{\psi_i\}_{i \in I_j}\entailsMod
                p$,
      that there are ($J$ need 
      not be finite). 
      \\
      multicut &  $p^n$ stands for $n$ occurrences of $p$; $m,n > 0$  
    \end{tabular}
    \rule{\linewidth}{\proofrulebreadth}
  \end{center}
\end{table}

\begin{remark}
  The following considerations may make the calculus, and the
  connection between $\langMod$ and $\lang$, more perspicuous. If we
  were to assume that $\lang$ has arbitrary disjunctions, our cut
  elimination result would show that we could write
  \begin{displaymath}
    \nec p \: = \: \bigvee_{\psi_1, \ldots, \psi_k \entailsMod p}
    \phi_i\land \ldots \land \phi_k
  \end{displaymath}
  In this case, the left and right
  rules for $\nec$ would simply be the left and right rules
  for a disjunction $\bigvee_{i \in I} P_i$,
  namely
  \begin{equation}
    \begin{prooftree}
      \left \{ \Gamma, p_i \entailsMod \Delta\right\}_{i \in I}
      \justifies
      \Gamma, \bigvee_{i \in I} p_i \entailsMod \Delta
      \using
      \bigvee \text{L}
    \end{prooftree}
    \quad \text{and}\quad
    \begin{prooftree}
      \Gamma \entailsMod p_i, \Delta \quad (i \in I)
      \justifies
      \Gamma \entailsMod \bigvee_{i \in I} p_i 
      \using 
      \bigvee \text{R}
    \end{prooftree}
  \end{equation}
  This explains the form of the side conditions $\psi_1, \ldots, \psi_k
  \entailsMod p$ in the left and right rules for $\nec X$: they
  take the place of $i \in I$ in the left and right rules
  for $\bigvee_{i \in I}$. Just as in those latter left and right
  rules, the conditions $\psi_1, \ldots, \psi_k
  \entailsMod p$ occur positively in both left and right rules: they
  are side conditions and not premises. 
\end{remark}

\begin{remark}
  We should remark that the rules for $\nec$ are, \emph{in general},
  infinitary (and even when finite, the set of premises of
  $\nec\text{L}$ can well be undecidable). The qualification ``in
  general'' is important here: cut elimination will show that, for
  typical applications, this rule very often has a finite, decidable
  (and, indeed, very tractable) set of premises.  Indeed, many
  applications -- for example, that described in
  Appendix~\ref{sec:turn-emphl-univ} -- will only use the right rule
  for $\nec$, which is much more tractable.  Furthermore, even when
  the rules are actually infinitary, the system is still very useful
  for metatheoretical purposes: after all, the proof-theoretic use of
  infinitary systems has a very long and respectable history (see, for
  example, \citend[p.~164]{fairtlough98:_hierar_provab_recur_funct}).
\end{remark}
\begin{example}
  Here is an example of how the left rule might plausibly be
  used. Michael Scriven 
  \cite[p.~61]{scriven88:_explan_predic_and_laws}
  describes the following argument pattern: we know that event $p$
  occurred. We know that the causes of $p$ are $a, b,$ or $c$: but we 
  also know that neither $a$ nor $b$ occurred. So we conclude $c$.

  In our system, this would go as follows. We have $p$, and we are supposing
  the world to be explanatorily closed, so we can conclude $\nec p$ (here
  we are doing causal explanation, so our explanations really are causes). 
  We can also conclude that 
  \begin{displaymath}
    \nec p \lequiv a \lor b \lor c;
  \end{displaymath}
  we also have $\lnot b$ and $\lnot c$, so we conclude $a$. 
\end{example}

\section{Consequences of Cut Elimination}
The proof of cut elimination is technical (it relies on 
ordinal induction), and is given in the appendix. It has
important consequences, however, and we give some of
them here.

\begin{corollary}\label{corollary:conservativity}
The modal theory is a conservative extension of the non-modal
theory: that is, if $\Gamma$ and $\Delta$ are non-modal,
we have 
\begin{displaymath}
  \Gamma \entails \Delta \quad\text{iff}\quad \Gamma \entailsMod \Delta. 
\end{displaymath}
\end{corollary}
\begin{proof}
Any proof of $\Gamma \entails \Delta$ is a proof of 
$\Gamma \entailsMod \Delta$; however, any cut free proof of
$\Gamma \entailsMod \Delta$ cannot involve the modal rules, and
is thus a proof of $\Gamma \entails \Delta$. 
\end{proof}
So, from now on we can safely ignore the distinction between
$\entails$ and $\entailsMod$.

\subsection{Semantics}\label{section:semantics}
We can now prove a soundness and completeness theorem for our
logic. We do this as follows: we first identify our model
$\kripke{\Theta}$ of Definition~\ref{def:canMod}
as the canonical model of this modal logic. 
That is, semantic validity in this single
model will entail syntactic validity.  We can, then, prove soundness
by showing that any syntactically valid sequent is also semantically
valid when interpreted in the canonical model, and prove completeness
by showing the converse. Recall that $\kripke{\Theta}$ has, as worlds,
models of the non-modal language $\lang$, and that the 
accessibility relation is given by the relation $\causRel{\Theta}$ 
of Definition~\ref{def:causRel}. The corresponding modelling
relation between worlds of $\kripke{\Theta}$ and sentences of $\langMod$
will, as in Definition~\ref{def:canMod}, be written $\modelsNec$.

For convenience of notation, we identify models of $\lang$ with the
sets of nonmodal sentences true in them, and these are just the
maximal $\entails$-consistent sets of nonmodal sentences: let
$\allModels$ be the set of all such models.  Clearly, for $A$ nonmodal
and $M \in \allModels$, we have
\begin{displaymath}
  M \models a\quad\text{iff} \quad 
  a \in M \quad \text{iff} \quad
  M \entails a. 
\end{displaymath}
For this section,
letters $a, a', \ldots$ will range over nonmodal sentences,
whereas $p, p', \ldots$ will range over the whole of $\langMod$.

The proof of the following is immediate:
\begin{lemma}
For any $M$ and any non-modal $a$,
    $M \modelsNec a$ iff $M \models a$.
\end{lemma}

We now define a semantic entailment relation, $\semEntMod$: this
will be like the usual notion of modal semantic entailment,
but with respect to the single model $\kripke{\Theta}$. As we have
explained, $\kripke{\Theta}$ will turn out to be the canonical
model of our logic, so the restriction to $\kripke{\Theta}$ is
harmless. 

We first define semantic values for sentences in $\langMod$.
\begin{definition}
  For sentences $p$ in $\langMod$, define their \emph{semantic values}
  $\semval{p} \subseteq \allModels$ as follows:
    \begin{align*}
      \semval{a} & \quad=\quad 
      \{ M | M \forces a\} & \text{for $a$ a nonmodal atom}\\
      \semval{\lnot p} &\quad=\quad \allModels - \semval{p} \\
      \semval{p\land q} &\quad=\quad \semval{p} \cap\semval{q} \\
      \semval{p\lor q} &\quad=\quad  \semval{p}  \cup \semval{q}\\
      \semval{\nec p} & \quad=\quad \{ M |
        \forall M'. M \causRel{\Theta} M' \Rightarrow 
        M' \in \semval{q}\}
    \end{align*}
\end{definition}
The following lemma is immediate:
\begin{lemma}
  For $a$ nonmodal, $\semval{a} = \{ M \in \allModels | M \models a\}$
\end{lemma}

\begin{definition}
  For $\Gamma, \Delta \subseteq \langMod$, we say that 
  \begin{align*}
    \Gamma \semEntMod \Delta \quad\text{iff}\quad
    \bigcap_{p \in \Gamma}\semval{p}
    \subseteq
    \bigcup_{q \in \Delta} \semval{q}.
  \end{align*}
\end{definition}
We also need to generalise our operation $(\cdot)^{\Theta}$ of \eqref{eq:causExpModel}
so that it can be applied to general subsets, rather
than simply models, of $\lang$. We recall that, for $M \in \allModels$,
$M^{\Theta} \:\:\eqdef\:\: \closure{\{\psi| 
  M \models \phi \:\text{for some}\: \phi \causes \psi\}}$.
\begin{definition}
  Let $S \subseteq \lang$: define 
  \begin{displaymath}
    S^{\Theta} \quad \eqdef \quad 
    \closure{
      \{
      a \in \lang \:|\: S \entails \phi_1\lor \ldots\lor \phi_k
      \:\: \text{for rules}  \:\:
      \phi_i \causes \psi_i,
      \:\: \text{with}\:\: 
      \psi_i \entails a 
      \:\text{for all}\: i
      \}
    }
  \end{displaymath}
\end{definition}
Note that, if $S$ is a model, the new definition coincides with the old one.
We need the disjunction of $\phi$s in order to prove the following lemma,
which will be crucial:
\begin{lemma}\label{lemma:intersection}
  If $\{S_i\}_{i \in I}$ is a set of subsets of $\lang$, with each
  $S_i$ closed under entailment, we have
  \begin{displaymath}
    \bigcap_{i \in I} (S_i^{\Theta}) \quad = \quad \left(
      \bigcap_{i \in I} S_i
    \right)^{\Theta}
  \end{displaymath}
\end{lemma}
\begin{proof}
  The right to left containment is trivial. For the left to right
  containment, suppose that $a \in S_i^{\Theta}$ for all $i$. Thus, for
  all $i$, we have 
  \begin{align*}
    S_i &\entails \phi_{i_1}, \ldots, \phi_{i_{k_i}} \quad
    \text{with}\quad \psi_{i_j} \entails a\: \:\text{for all} \:\:j.
    \intertext{Then, since the $S_i$ are closed under entailment, we
      have}
    \bigcap_{i \in I}S_i & \entails \bigcup_{i \in I}\{\phi_{i_1},
    \ldots, \phi_{i_{k_i}}\}.
    \intertext{But, by the compactness of $\lang$, we can find a finite
      set of $\{\phi_i \causes \psi_i \}$ with}
    \bigcap_{i \in I} S_i & \entails \phi_1, \ldots, \phi_n \quad 
    \text{with}\quad  \psi_{i} \entails a\: \:\text{for all} \:\:i,
  \end{align*}
  and this establishes the result. 
\end{proof}

We now prove soundness and completeness.

\begin{lemma}\label{lemma:restrictedSoundComplete}
  The $\kripke{\Theta}$ semantics is sound and complete for 
  inferences of the form $S \entailsMod \nec p$, for
  $\Gamma \subseteq \lang$ and $p \in \langMod$: that is,
  \begin{displaymath}
    \Gamma \entailsMod p \quad \text{iff}\quad \Gamma \semEntMod p.
  \end{displaymath}
\end{lemma}
\begin{proof}
  We prove this by induction on the complexity of $p$: for
  $p$ nonmodal, it is, by Corollary~\ref{corollary:conservativity},
  simply soundness and completeness for the predicate calculus. 
  The inductive steps for the nonmodal connectives are quite standard:
  there remains the case when $p = \nec q$. It will suffice
  to show that 
  \begin{align*}
     \semval{\nec q}\quad & = \quad\bigcup_i \semval{\phi_{i, 1} \land
      \ldots \land \phi_{i, k_i}}
    \intertext{where, in the usual way, $\phi_{i,j} \causes \psi_{i,j}$
      for all relevant $i,j$, where $\psi_{i,1} , \ldots , 
      \psi_{i,k_i} \entailsMod q$ for all $i$, and where these are
      all of such finite sets of causal laws that there are; the
      right to left inference gives soundness, and the left
      to right inference gives completeness. Note that the semantic values of
  $\nec q$, for any $q$, are a union of the semantic values of nonmodals, and so}
    \semval{\nec q} & = \bigcup_{\semval a \subseteq \semval{\nec q}} \semval{a}
      \intertext{Now, if $a$
        is nonmodal,}
    \semval{a}\quad &\subseteq\quad \semval {\nec q} \\
    & \text{iff}\quad
    \forall M. M \models a . \forall M'. M' \forces M^{\Theta}
    \:\:\Rightarrow \:\: M' \modelsNec q
    \intertext{expanding the definitions of $\semval{\nec q}$
      and of $\causRel{\Theta}$}
    &\text{iff}\quad \forall M'. (\exists M \models q. M' \forces M^{\Theta})
    \:\:\Rightarrow\:\: M' \modelsNec q\\
    & \text{iff}\quad \forall M'. M' \forces \left(\closure{\{A\}}\right)^{\Theta}
    \:\:\Rightarrow \:\: M' \modelsNec q
    \intertext{(where $\closure{\{a\}}$ is the deductive closure of $a$)
      by Lemma~\ref{lemma:intersection}}
    &\text{iff}\quad \left(\closure{\{a\}}\right)^{\Theta}
    \modelsNec q
    \intertext{by the inductive hypothesis. Now expansion of
      the definitions shows that}
    \bigcup_{\left(\closure{\{a\}}\right)^{\Theta}\modelsNec Q} 
    \semval{a} \quad &
    = \quad\bigcup_i \semval{\phi_{i, 1} \land
      \ldots \land \phi_{i, k_i}}
  \end{align*}
  (with the usual restrictions on the $\phi_{i,j}$), which was to be proved. 
\end{proof}
Recall now that the forcing relation, $\modelsNec$, on the Kripke 
model is defined inductively using the structure of $\kripke{\Theta}$. 
We can now show that, in fact, it is the same as provability. 
\begin{corollary}\label{corollary:redefineForcing}
  For a model $M \in \allModels$ and $p \in \langMod$,  
  \begin{displaymath}
    M \modelsNec p \quad \text{iff} \quad M \entailsMod p.
  \end{displaymath}
\end{corollary}
\begin{proof}
  By soundness and completeness, 
  \begin{align*}
    M \entailsMod p \quad&\text{iff} \quad \semval{M} \subseteq \semval{p};
    \intertext{but $\semval{M} = \{ M \}$, so}
    M \entailsMod p \quad&\text{iff} \quad M \in \semval {p}\\
    &\text{iff}\quad M \modelsNec p
  \end{align*}
  by definition of $\modelsNec$
\end{proof}
\begin{theorem}
  The $\kripke{\Theta}$ semantics is sound for general sequents.
\end{theorem}
\begin{proof}
  We prove this by induction over the length of a proof: the nonmodal
  rules are straightforward, and we use
  Lemma~\ref{lemma:restrictedSoundComplete} for the modal rules.
\end{proof}
\begin{lemma}\label{lemma:defExt}
  For any $p \in \langMod$, and for any $M \in \allModels$, 
  either $M \entailsMod p$ or $M \entailsMod \lnot p$.
\end{lemma}
\begin{proof}
  We prove this by induction: since $\entailsMod$ coincides
  with $\entails$ on $\lang$, and since $M$ is a model of $\lang$,
  the result is clearly true if $p$ is nonmodal. The crucial step
  is to prove the result for propositions of the form $\nec q$: 
  here we argue as follows
  \begin{align*}
    \text{Suppose} && M &\not\entailsMod \lnot \nec q\\
    \text{then}&& M, \nec q \quad&\not\entailsMod\quad \lfalse
    \intertext{and, by the left rule for  $\nec$, 
      this holds iff, for some set of causal rules $I$,}
    &&M, \{\phi_i\}_{i \in I} \quad&\not\entailsMod\quad \lfalse
    \intertext{which, since $M$ is a model of $\lang$, in turn holds iff,
      for any $i\in I$,} 
    &&M \quad&\entailsMod\quad \phi_i&
    \intertext{and so, by the right rule for $\nec$,}
    &&M \quad &\entailsMod\quad \nec q 
  \end{align*}
  By cut elimination, no element of $\langMod$ can
  be both entailed and not entailed by $M$. 
\end{proof}

\begin{corollary}
  The maximal $\entailsMod$-consistent subsets of $\langMod$
  are in one-to-one correspondence with the maximal $\entails$-consistent
  subsets of $\lang$ (i.e.\ the elements of $\allModels$).
\end{corollary}
\begin{proof}
  Intersection with $\lang$ gives a morphism in one direction:
  Lemma~\ref{lemma:defExt} shows that every maximal $\entails$-consistent
  subset of $\lang$ can be uniquely extended to a maximal
  $\entailsMod$-consistent subset of $\langMod$. 
\end{proof}
\begin{theorem}
  The $\kripke{\Theta}$ semantics is complete.
\end{theorem}
\begin{proof}
  It suffices to show that, if $\Gamma \subseteq \langMod$, and if
  $\Gamma \not\entailsMod \lfalse$, then $\semval{\Gamma}$ is not
  empty, i.e.\ that there is some world $M$ of $\kripke{\Theta}$ with $M
  \in \semval{\Gamma}$.  But, if we have $\Gamma \not\entailsMod
  \lfalse$, we can extend $\Gamma$ to a maximal
  $\entailsMod$-consistent subset of $\langMod$: this will correspond
  to an element $M$ of $\allModels$. By
  Corollary~\ref{corollary:redefineForcing}, this is the world we
  want.
\end{proof}
\subsubsection{Interpolation and Finiteness}
We first prove the following interpolation theorem. 
\begin{proposition}[Interpolation]\label{thm:interpolation}
  If we have $\Gamma \entailsMod \nec p, \Delta$,
  for some $p$ and $\Gamma, \Delta \subseteq \langMod$,
  then there is a set $\Delta_0 \subseteq \lang$ 
  such that
  \begin{align*}
    \Gamma &\entailsMod \Delta_0, \Delta
    \intertext{and, for each $a \in \Delta_0$, there is some $b \in \lang$
      such that}
    a & \entailsMod \nec b&\text{and}\\
    b &\entailsMod p. 
  \end{align*}
\end{proposition}
\begin{proof}
  We prove this by a straightforward induction over a 
  cut free proof: because we might use the contraction rule
  on $\nec p$, the inductive hypothesis will be the same
  as that of the theorem but with the slightly weaker 
  condition that $\Gamma \entailsMod (\nec p)^n, \Delta$.
  The elements $a$ of  $\Delta_0$ will turn out to be finite
  conjunctions of rule bodies, whereas the $b$s
  will be the corresponding conjunctions of rule heads. 
\end{proof}
One of the key applications of our cut elimination theorem
is in proving finiteness results, such as the following: it is
applicable to the modal implications that we will use 
in the next section for interpreting Parsons and Jennings' 
logic of argument, so it is quite significant. 
\begin{corollary}\label{cor:normForm}
  If 
  \begin{displaymath}
    \Gamma \entails \nec p
  \end{displaymath}
  is provable in the sequent calculus, with $\Gamma$ a set of nonmodals, then
  there is a finite set $\Gamma' \subseteq \Gamma$, and nonmodal
  propositions $A$ and $B$, such that 
  \begin{align*}
    \Gamma' & \entails a,\\
    a &\entailsMod \nec b,&\text{and}\\
    b &\entailsMod p.
  \end{align*}
\end{corollary}
\begin{proof}
  We apply Theorem~\ref{thm:interpolation}, which gets us a possibly
  infinite set $\Delta$ with $\Gamma \entailsMod\Delta$. But, in this
  case, $\Gamma$ and $\Delta$ are both classical, so we can apply
  Corollary~\ref{corollary:conservativity} to show that $\Gamma
  \entails \Delta$. We can now apply compactness to find a \emph{finite}
  $\Delta_0 \subseteq \Delta$, and, taking disjunctions, we get a single
  $a$ and $b$. 
\end{proof}

\section{Explanation and Argument}
As we have argued, McCain and Turner's theory seems to be a very
natural formulation of explanation \emph{in general}: although its
original application may have been to a causal context, there is
nothing about it which forces these explanations to be
\emph{causal}. Once we broaden our horizons to general explanation, we
can bring this modal system into contact with other work: we give one
example of this, the theory of Parsons and Jennings
\cite{parsons:_negot_argum} (see also \citend{parsons98agents}). We
could give other applications: for example, where we have
explanations, questions must also be in the neighbourhood, and thus we
can also express a good deal of the formalism of
\citend{belnap76:_logic_quest_answer} (see \citend{harrah98:_quest})
in terms of ours.

\subsection{The Parsons and Jennings System}
Parsons and Jennings (\cite{parsons:_negot_argum}; see also
\citend{parsons98agents}) have described a
consequence relation, $\entailsACR$, which is
intended to capture the practice of argumentation.
The items which this system manipulates are ordered pairs
$(p,\Gamma)$, where $p$ is a proposition and $A$ a
set of propositions: we will call such a 
pair an \emph{argument}. Intuitively, $p$ is the conclusion of
an argument, and $\Gamma$ is the set of its grounds. We will also
call $p$ the \emph{head} of the  argument, and $\Gamma$ its
\emph{body}.

The system is given in Table~\ref{tab:parsonsJennings}: here
$\Theta$ is a set of \emph{basic} arguments, and it plays exactly the
same role in their system as a causal theory does in McCain and Turner's.
Parsons and Jennings write their system 
in natural deduction style, 
with introduction and elimination rules, using 
sequents of the form
$\Theta \entailsACR (p, A)$: such a sequent says that
$A$ is an argument for $p$, given basic arguments $\Theta$.

(Note that we have interchanged the labels on the rules
$\lnot \text{I}$ and RAA: Parsons and Jennings' original labelling
is clearly a typo of some sort.)  
\begin{table}[ht]\centering
  \setlength{\extrarowheight}{23pt}
  \caption{The Parsons and Jennings System}
  \label{tab:parsonsJennings}
  \rule{\linewidth}{\proofrulebreadth}\\
  \begin{displaymath}
    \begin{array}{cc}
      \begin{prooftree}
        (p, \Gamma) \in \Theta
        \justifies
        \Theta \entailsACR (p, \Gamma)
        \using
        \text{Axiom}
      \end{prooftree}
      &
      \begin{prooftree}
        \vphantom{(p, \Gamma) \in \Theta}
        \justifies
        \Theta \entailsACR (\ltrue, \emptyset)
        \using
        \ltrue \text{I}
      \end{prooftree}
      \\
      \begin{prooftree}
        \Theta \entailsACR (p,\Gamma)
        \quad
        \Theta \entailsACR (q,\Gamma')
        \justifies
        \Theta \entailsACR (p \land q, \Gamma \cup \Gamma')
        \using
        \land\text{I}
      \end{prooftree}
      &
      \begin{array}{c}
        \begin{prooftree}
          \Theta \entailsACR (p \land q, \Gamma)
          \justifies
          \Theta \entailsACR (p,\Gamma)
          \using
          \land \text{E1}
        \end{prooftree}
        \\[-9pt]
        \begin{prooftree}
          \Theta \entailsACR (p \land q, \Gamma)
          \justifies
          \Theta \entailsACR (q,\Gamma)
          \using
          \land \text{E2}
        \end{prooftree}
      \end{array}
      \\
      \begin{array}{c}
        \begin{prooftree}
          \Theta \entailsACR (p,\Gamma)
          \justifies
          \Theta \entailsACR (p \lor q, \Gamma)
          \using
          \lor \text{I1}
        \end{prooftree}
        \\[-9pt]
        \begin{prooftree}
          \Theta \entailsACR (q,\Gamma)
          \justifies
          \Theta \entailsACR (p \lor q, \Gamma)
          \using
          \lor \text{I2}
        \end{prooftree}
      \end{array}
      &
      \begin{array}{c}
        \begin{prooftree}
          {\setlength{\extrarowheight}{2pt}
            \begin{array}{l}
              \Theta\entailsACR (p \lor q, \Gamma) \\
              \quad\Theta, (p,\Gamma) \entailsACR (r, \Gamma')\\
              \qquad\Theta, (q,\Gamma) \entailsACR (r,\Gamma'') 
            \end{array}
          }
          \justifies
          \Theta \entailsACR (r, \Gamma' \cup \Gamma'')
          \using
          \lor \text{E}
        \end{prooftree}
      \end{array}
      \\
      \begin{prooftree}
        \Theta, (p, \emptyset) \entailsACR (\lfalse, \Gamma)
        \justifies
        \Theta \entailsACR (\lnot p,\Gamma)
        \using
        \lnot \text{I}
      \end{prooftree}
      &
      \begin{prooftree}
        \Theta \entailsACR (p, \Gamma) 
        \:
        \Theta \entailsACR (\lnot p, \Gamma)
        \justifies
        \Theta \entailsACR (\lfalse, \Gamma)
        \using
        \lnot \text{E}
      \end{prooftree}
      \\
      \begin{prooftree}
        \Theta, (p, \emptyset) \entailsACR (q,\Gamma)
        \justifies
        \Theta \entailsACR(p \rightarrow q,\Gamma)
        \using
        \rightarrow \text{I}
      \end{prooftree}
      &
      \begin{prooftree}
        \Theta \entailsACR (p, \Gamma) 
        \:
        \Theta \entailsACR (p \rightarrow q, \Gamma)
        \justifies
        \Theta \entailsACR (q,\Gamma \cup \Delta)
        \using
        \rightarrow \text{E}
      \end{prooftree}
      \\
      \begin{prooftree}
        \Theta \entailsACR (\lfalse, \Gamma)
        \justifies
        \Theta \entailsACR (p, \Gamma)
        \using
        \text{EFQ}
      \end{prooftree}
      &
      \begin{prooftree}
        \Theta, (\lnot p,\emptyset) \entailsACR (\lfalse, \Gamma)
        \justifies
        \Theta \entailsACR (p, \Gamma)
        \using
        \text{RAA}
      \end{prooftree}
    \end{array}
  \end{displaymath}
  \rule{\linewidth}{\proofrulebreadth}
\end{table}

\subsection{Comparison with Our System}
We can now translate Parsons and Jennings' system into ours. 
\begin{definition}
  Let $\Theta \entailsACR (p,\Gamma)$ be a sequent in Parsons and Jennings' system. 
  Its \emph{modal translation} is the sequent 
  \begin{displaymath}
    \Gamma \entails \nec_{\Theta} P,
  \end{displaymath}
  where $\nec_{\Theta}$ is the modal operator defined by causal laws 
  \begin{displaymath}
    \left\{
      \left(\bigwedge_{q \in \Gamma}q\right) \causes p \bigm| (p,\Gamma) \in \Theta
    \right\}.
  \end{displaymath}
\end{definition}

Since the Parsons and Jennings system is written in natural deduction
style, some of the rules (for example $\rightarrow\text{I}$) manipulate
the set of basic arguments: consequently, the modal operator in the 
modal translation will vary. We will, then, need the following
lemma: 
\begin{lemma}\label{lem:modalityExpansion}
  If $\nec$ is the modality associated to a set $\Theta$ of
  explanations, and if $\nec'$ is
  that associated to $\Theta \cup \{(\psi,\phi)\}$, then define
  an interpretation 
  \begin{align*}
  \alpha: \lang_{\nec'} &\rightarrow  \lang_{\nec}
  \intertext{by}
  \alpha(p) &= p & \text{for $p$ atomic}\\
  \alpha(p\land p') & = \alpha(p) \land \alpha(p')\\
  \alpha(p \lor p' ) & = \alpha(p) \lor \alpha(p')\\
  \alpha(\lnot p) & = \lnot \alpha(p) \\
  \alpha(\nec' p) & = 
    (\phi \land \nec
    (\psi \rightarrow \alpha(p))) \lor \nec \alpha(p)
  \end{align*}
  for any  $p$. Then, for any $\Gamma, \Delta \subseteq \lang_{\nec'}$
  \begin{displaymath}
    \Gamma \entails_{\nec'} \Delta \quad\text{iff} 
    \quad 
    \alpha(\Gamma) \entailsMod \alpha(\Delta)
  \end{displaymath}
\end{lemma}
\begin{proof}
  We check that $\Gamma$ and $\alpha(\Gamma)$ have the
  same semantic values (regarded as subsets of $\allModels$),
  and similarly for $\Delta$; we then apply the soundness and completeness
  theorems. 
\end{proof}

We have the following
\begin{proposition}\label{proposition:modTransSound}
  The modal interpretation is \emph{sound}: that is, each
  of Parsons and Jennings' axioms is translated into
  a tautology. 
\end{proposition}
\begin{proof}
  \begin{description}
  \item[Ax] is
    \begin{displaymath}
      \begin{prooftree}
        \justifies
        \Theta \entailsACR (p,\Gamma)
        \using
        (p,\Gamma) \in \Theta
      \end{prooftree}
    \end{displaymath}
    and this follows from our definition of the modal translation.
    
  \item[$\land$--I, $\lnot$--E, $\rightarrow$--E]
    $\land-\text{I}$, for example,  is 
    \begin{displaymath}
      \begin{prooftree}
        \Theta \entailsACR (p,\Gamma) 
        \quad 
        \Theta \entailsACR (q,\Gamma')
        \justifies
        \Theta \entailsACR (p \land q, \Gamma \cup \Gamma')
      \end{prooftree}
    \end{displaymath}
    and this
    follows from the 
    \textbf{K} tautology $\nec p \land \nec q \entails \nec p \land q$. 
    $\lnot-\text{E}$ and  $\rightarrow-\text{E}$  are similar. 
    
  \item[$\land$--E1, $\land$--E2, $\lor$--I1, $\lor$--I2, EFQ]
    $\land-\text{E1}$, for example, is 
    \begin{displaymath}
      \begin{prooftree}
        \Theta \entailsACR (p \land q,\Gamma )
        \justifies
        \Theta \entailsACR (p,\Gamma)
      \end{prooftree}
    \end{displaymath}
    and this follows from the \textbf{K} tautology
    $\entails \nec(a \land b) \rightarrow  \nec a$.
    $\land-\text{E2}$, $\lor-\text{I1}$,  $\lor\text{I2}$, and EFQ are similar.
    
  \item[$\ltrue$--I]This is just
    \begin{displaymath}
      \entails \nec \ltrue,
    \end{displaymath} 
    a \textbf{K} tautology. 
    
  \item[$\lor$--E]This is
    \begin{displaymath}
      \begin{prooftree}
        \Theta \entailsACR (p \lor q, \Gamma)
        \quad
        \Theta, (p,\Gamma) \entailsACR (r,\Gamma')
        \quad
        \Theta, (q,\Gamma) \entailsACR (r,\Gamma'')
        \justifies
        \Theta \entailsACR (r, \Gamma' \cup \Gamma'') 
      \end{prooftree}
    \end{displaymath}
    and this corresponds, in our system, to
    \begin{multline*}
      \qquad \qquad
      a \rightarrow \nec_{\Theta} (p \lor q),\,
      b \rightarrow \nec_{\Theta \cup \{(p,\Gamma)\}}r,\,
      c \rightarrow \nec_{\Theta \cup \{(q,\Gamma)\}}r
      \\
      \entails
      b \land c \rightarrow \nec_{\Theta} r
    \end{multline*}
    We can use the lemma to express $\nec_{\Theta \cup\{(p,\Gamma)\}}$ and 
    $\nec_{\Theta \cup \{(q,\Gamma)\}}$ in terms of $\nec_{\Theta}$; some
    routine but tedious computation then reduces this case to
    \begin{equation}
      \nec_{\Theta} (p \lor q), \nec_{\Theta} (p \rightarrow r ), 
      \nec_{\Theta} (q \rightarrow r)
      \entails \nec_{\Theta} r
    \end{equation}
    which is a \textbf{K} tautology. 
    
  \item[RAA, $\lnot$--I]$\lnot\text{--I}$ is
    \begin{displaymath}
      \begin{prooftree}
        \Theta, (p,\emptyset) \entailsACR (\lfalse,A)
        \justifies
        \Theta \entailsACR (\lnot p,A)
      \end{prooftree}
    \end{displaymath}
    which corresponds to 
    \begin{displaymath}
      A \rightarrow \nec_{\Theta \cup \{(p,\emptyset)\}} \lfalse
      \entails
      A \rightarrow \nec_{\Theta} \lnot p;
    \end{displaymath}
    using the lemma on $\nec_{\Theta \cup \{(p,\emptyset)\}}$,
    and some computation, reduces this to
    \begin{displaymath}
      \nec \lfalse \lor \nec (p \rightarrow \lfalse) \entails \nec \lnot p
    \end{displaymath}
    which is a \textbf{K} tautology. $\lnot$--I is similar. 
    
  \item[$\rightarrow$--I]This  is
    \begin{displaymath}
      \begin{prooftree}
        \Theta, (p,\emptyset) \entailsACR (q,A)
        \justifies
        \Theta \entailsACR (p \rightarrow q,A)
      \end{prooftree}
    \end{displaymath}
    which corresponds to
    \begin{displaymath}
      A \rightarrow \nec_{\Theta \cup \{(p,\emptyset)\}} q
      \entails
      A \rightarrow \nec_{\Theta} (p \rightarrow q).
    \end{displaymath}
    The usual moves reduce this to 
    \begin{displaymath}
      \nec q \lor \nec (p \rightarrow q) \entails \nec (p \rightarrow q)
    \end{displaymath}
    again a \textbf{K} tautology. 
  \end{description}
\end{proof}

Completeness does not hold. This is for trivial reasons: all rules
(except $\lor\text{E}$) of the Parsons and Jennings system leave the
body of the argument intact. A trivial induction on the length of
proofs will yield
\begin{proposition}
  In any proof of 
  \begin{displaymath}
    \Theta \entailsACR (p,\Gamma),
  \end{displaymath}
  $A$ must be a union of the bodies of rules in $\Theta$. 
\end{proposition}
Since the modal sequent calculus certainly does not satisfy this
condition, we cannot hope for completeness. What we need to do
is to be able to compose proofs in the Parsons and Jennings system
with natural deduction proofs for the grounds of an argument: 
we could, theoretically, write down another set of rules for doing this. 
However, we only need one extra rule, which is this:

\begin{definition}
  Let \emph{classical $\lor\text{E}$} be the following rule:
  \begin{displaymath}
    \begin{prooftree}
      \Gamma \entails q \lor r
      \quad
      \Theta \entailsACR (r, \Gamma' \cup \{q \})
      \quad
      \Theta \entailsACR (r, \Gamma'' \cup \{ r\})
      \justifies
      \Theta \entailsACR (r, \Gamma \cup \Gamma' \cup \Gamma'')
      \using
      \lor\text{EC}
    \end{prooftree}
  \end{displaymath}
  where $A \entails q \lor r$ is an entailment in classical natural deduction. 
\end{definition} 

We clearly have
\begin{proposition}
  The modal translation is sound for classical or-elimination. 
\end{proposition}

And we can also prove completeness:
\begin{theorem}
  The modal translation is complete: that is, given a proof of
  \begin{equation}
    \Gamma \entails \nec_{\Theta} p,
  \end{equation}
  there is a proof, in the Parsons and Jennings system together with
  classical or-elimination, of
  \begin{equation}
    \Theta \entailsACR(p,\Gamma)
  \end{equation}
\end{theorem}

\begin{proof}[Sketch of proof]
  We establish the following  lemma:
  \begin{lemma}\label{lem:mirrorSeqPrf}
    If $p$ is non-modal, given a sequent calculus proof of
    \begin{displaymath}
      \psi_1, \ldots, \psi_k \entails p,
    \end{displaymath}
    then there is a Parsons and Jennings proof of
    \begin{displaymath}
      \Theta \entails (p, \{\phi_1, \ldots, \phi_k\}).
    \end{displaymath}
  \end{lemma}
  This lemma can be proved by first 
  transforming the sequent calculus proof to a natural deduction
  proof, and then observing that the Parsons and Jennings rules
  mirror the rules of classical natural deduction.
  
  So now we can prove the theorem: we take a proof of $A \entails
  \nec_{\Theta} p$, and from it derive a finite set $I$ such that
  \begin{align}
    \Gamma &\entails \bigvee_{i \in I} \phi_{i_1} 
    \land \ldots \land \phi_{i_{k_i}}\label{eq:2}\\
    \Theta &\entailsACR (p, \{ \phi_{i_1}, \ldots, 
    \phi_{i_{k_i}}\})&&\text{for any $i$}\label{eq:3}
  \end{align}
  We then use classical or-elimination in order to glue together
  \eqref{eq:2} and \eqref{eq:3}.
\end{proof}

\begin{remark}
  As we see here, the natural deduction formulation is actually quite
  ambiguous as to what its premises are: in a proof of
  \begin{displaymath}
    \Theta \entailsACR (p, \Gamma),
  \end{displaymath}
  are the premises the basic arguments $\Theta$, or the grounds $\Gamma$ for
  the argument which is to be established? Now the system is set up as
  if the premises are the set $\Theta$ of basic arguments, and this gives
  a sense of composition of arguments which is valid: that is, from
  $\Theta \entailsACR (p,\Gamma)$ and $\Theta', (p,\Gamma) \entailsACR (q,\Gamma')$, we can
  establish $\Theta', \Theta \entailsACR (q,\Gamma')$. But this is not enough: we
  \emph{also} want to regard the grounds of arguments as premises.
  
  The situation is clearly two-dimensional in something like Pratt's
  sense -- he defines the \emph{dimension} of a logic to be ``the
  smallest number of variables and constants of the logic sufficient
  to determine the remaining variables and constants''
  \citend{pratt93:_roadm_some_two_dimen_logic}: the modal operator can
  be varied quite independently of the classical connectives, merely
  by altering the set of causal rules. Consequently, a formalism such
  as Masini's \citend{masini92:_sequen_calculPT}
  \citend{masini93:_sequen_calculINT} may well be more appropriate.
\end{remark}

\begin{remark}
  This translation between sequent calculus and the Parsons and
  Jennings natural deduction is, in addition, not very sensitive to
  the structure of proofs on either side: natural deduction proofs
  tend to transform the conclusion of the argument quite extensively
  before coming down to basic arguments. Sequent calculus proofs, by
  contrast, leave the conclusion unchanged until an application of
  $\nec R$, after which the proof is a matter of standard classical
  logic.  In addition, the Parsons and Jennings system only represents
  a fragment of the full sequent calculus (namely, the entailments
  in which $\nec$ only occurs on the right). A natural deduction formulation
  of the entire sequent calculus would be interesting, but would
  have to extend the Parsons and Jennings system quite considerably. 
\end{remark}
\section{Conclusion: Two Approaches to Nonmonotonic 
  Reasoning}
There are two approaches to nonmonotonic reasoning: there is 
the generally accepted one, which may well be susceptible to
Fodor's critique. It can be summarised in Brewka's words:
\begin{quote}
  To formalise human commonsense reasoning something different [from
  classical logic] is needed. Commonsense reasoning is frequently not
  monotonic. In many situations we draw conclusions which are given up
  in the light of further information. \cite[p.~2]{brewka}
\end{quote}
According to this view, nonmonotonic logic is applicable \emph{globally}: 
it applies to \emph{all} of our commonsense reasoning. There is,
so to speak, a small fragment of our reasoning which happens to be
monotonic: but commonsense is nonmonotonic by default. There would, then,
be a large amount of reasoning which was both nonmonotonic and global. 
This, if true, would be fairly catastrophic. There has been quite a lot
of recent success in implementing nonmonotonic logic, but it is still
true that, precisely because of this globality, performance scales quite
badly \cite{cadoli94is}. So, if we are supposed to do it on a very
large scale, it cannot be expected to perform well. 

Contrast this with the McCain-Turner system. There are two aspects
that we need to examine: firstly, reasoning within the system itself,
and, secondly, the relation between the set of causal rules and the
modal operator.  The system itself is monotonic. It is a modal logic
of a fairly well-known sort, and it has a sequent calculus which
admits cut elimination; we can, then, search for proofs
efficiently. The dependence of the modal operator on the set of rules
$\phi \causes \psi$ is, however, non-monotonic: if we add new rules,
thus changing the modal operator, then we may invalidate modal
inferences that we previously made. However, only the \emph{left} rule
for $\nec$ is thus sensitive; the right rule is not. Consequently,
although our system \emph{is} nonmonotonic, it is so in a quite
limited way: it is not susceptible to the sort of ``anything might
entail anything else'' worries that Fodor has.

There is a final remark to be made. Many of the tractability 
results for our system do not follow directly from its definition, but
are established on the basis of metatheory, and quite technical metatheory
at that. Now this would, perhaps, be a fatal objection if we were 
supposed to adopt this system by means of introspection. However, 
if we are supposed to acquire styles of reasoning on the basis of
natural selection -- in the manner of the cognitive science 
described in \cite[ch.~5]{FodorJA:mindww} -- then none of this matters:
the metatheoretical results apply to the system, and give it the
evolutionary advantages that it arguably has, regardless of
whether or not the early primates who, maybe, adopted these styles
of reasoning were, or were not, familiar with technical results in
proof theory. 

\newpage
\appendix

In the appendix, we have material that does not fit naturally 
into the argumentative structure of the main paper: the 
cut elimination, which is too technical, and some material on
Turner's logic of universal causation, which relates what we have
done to previous work. 
\section{Cut Elimination}
We now prove our cut elimination result. This needs a certain amount
of machinery, because our system is, in generally, infinitary: proof
trees may, therefore, be infinitely branching, and -- since we cannot
apply Zorn's lemma -- there may be branches of infinite length above a
given node. Because of this, we cannot assign a finite depth to each
node, and, consequently, we cannot prove cut elimination in the usual
manner, that is, by an induction on both the depth and the complexity
of the cut formula.  We can use an induction, but it must be an
induction over ordinals: there is a machinery of ordinal analysis,
which is used for the proof-theoretic analysis of infinitary systems
(for example, arithmetic with the $\omega$-rule). We will generally
follow the treatment in \citend{pohlers89:_proof_theor}, with the
modifications necessary because our system is two-sided, not
one-sided, and because our sole infinitary connective is $\nec$ rather
than the infinite conjunctions and disjunctions which Pohlers
investigates. Our case, in fact, is rather simple, because, if we
regard $\nec$ as a possibly infinite disjunction of finite conjuncts of
non-modal formulae, the complexity of the disjuncts
is bounded \emph{a priori} by $\omega$.

In this section, since we want to treat entailments uniformly, we will
write rule applications as follows:
\begin{displaymath}
  \begin{prooftree}
    \left\{\Gamma_i \entailsMod \Delta_i\right \}_{i \in I}
    \justifies
    \Gamma \entailsMod \Delta
  \end{prooftree}
\end{displaymath} 
where $I$ can have zero, one, or two members (for the
standard rules) or an arbitrary number (for $\nec\text{L}$).

We first define the \emph{rank} of a formula: 
\begin{definition}
  If $F \in \langMod$, let its rank $\rank(F)$ be the
  ordinal defined inductively as follows:
  \begin{enumerate}
  \item $\rank(F) = 0$ if $F$ is atomic, 
  \item $\rank (\lnot F) = \rank(\forall x. F) = \rank (\exists x. F) 
    = \rank(F) + 1$
  \item $\rank (F_1 \land  F_2) = \rank (F_1 \lor F_2) =
    \rank (F_1 \rightarrow F_1) = \max (\rank(F_1), \rank(F_2)) + 1$
  \item $\rank (\nec F) = 0.$
  \end{enumerate}
\end{definition} 
Note that $\rank (\nec F)$ is independent of $F$: this is
because, as we have remarked, we can give an \emph{a priori} bound
on the disjuncts that $\nec F$ can be thought of as
expanding into. We have
\begin{proposition}\label{proposition:immediateRankBound}
  For any $F \in \langMod$, $\rank(F) < \omega$. 
\end{proposition}
\begin{proof}
Immediate. 
\end{proof}

We define the following classification of inference rules:
\begin{definition}
  The \emph{finitary} rules will be all of the rules apart from
  $\nec\text{L}$. The \emph{non-cut} rules will be
  all of the rules except cut.
\end{definition}
Next, we define entailment symbols annotated with both the rank of
cutformulae and the depth of the proof tree. There are two classes of
inference rules: the finitary rules, and $\nec\text{L}$.
$\nec\text{L}$ is the only infinitary rule; our proof trees,
then, will consist of portions where finitary rules are applied,
separated by applications of $\nec\text{L}$. We will treat
the two kinds of rule application differently, and
we will, therefore, use
two ordinals for measuring the depth of the tree. 
$\zeta$ will measure
the number of applications of $\nec\text{L}$ above the root of the
tree which have proofs involving cuts above them: it 
will, in general, be infinite. $\alpha$ will measure the
number of finitary rule applications above the root of the tree, but
will be reset to zero after each application of $\nec\text{L}$: it
will be finite. $\rho$ will measure the rank of cut formulae, and will
likewise be reset to zero after any application of $\nec\text{L}$:
it will, therefore, be finite. The
last two rules will make the annotations fit smoothly with our
inductive arguments, which will first over $\zeta$, then over
$\rho$, then over $\alpha$. 
\begin{definition}\label{definition:entailsRanked}
  Define the entailment relations
  $\entailsRanked{\alpha}{\zeta}{\rho}$ as follows:
  \begin{enumerate}
  \item If $\Gamma \entailsMod \Delta$ is an instance of $\text{Ax}$, of
    $\lfalse\text{L}$, or of $\ltrue\text{R}$, then
    \begin{displaymath}
      \Gamma \entailsRanked{0}{0}{0} \Delta
    \end{displaymath}
  \item If we have an instance of a finite non-cut 
    inference rule
    \begin{displaymath}
      \begin{prooftree}
        \left\{\Gamma_i \entailsMod \Delta_i\right\}_{i \in I}
        \justifies
        \Gamma \entailsMod \Delta
      \end{prooftree}
    \end{displaymath}
    and if, for all $i \in I$, 
    $\Gamma_i \entailsRanked{\alpha_i}{\zeta}{\rho} \Delta_i$
    where, for all $i$, 
    $\alpha_i < \alpha$, then 
    \begin{displaymath}
      \Gamma \entailsRanked{\alpha}{\zeta}{\rho} \Delta. 
    \end{displaymath}
  \item If 
    \begin{displaymath}
      \begin{prooftree}
        \left\{\Gamma_i \entailsMod \Delta_i\right\}_{i \in I}        
        \justifies
        \Gamma \entailsMod \Delta
      \end{prooftree}
    \end{displaymath}
    is an instance of $\nec\text{L}$, and if, for all $i$,
    $\Gamma_i \entailsRanked{\alpha_i}{0}{0} \Delta_i$
    where $\zeta_i < \zeta$, then
    \begin{displaymath}
      \Gamma \entailsRanked{0}{0}{0} \Delta.
    \end{displaymath}
  \item If 
    \begin{displaymath}
      \begin{prooftree}
        \left\{\Gamma_i \entailsMod \Delta_i\right\}_{i \in I}        
        \justifies
        \Gamma \entailsMod \Delta
      \end{prooftree}
    \end{displaymath}
    is an instance of $\nec\text{L}$, and if, for all $i$, $\Gamma_i
    \entailsRanked{\alpha_i}{\zeta_i}{\rho_i} \Delta_i$ where either
    some $\zeta_i > 0$ or some $\rho_i > 0$ and where $\zeta_i <
    \zeta$ for all $i$, then
    \begin{displaymath}
      \Gamma \entailsRanked{0}{\zeta}{0} \Delta.
    \end{displaymath}
  \item If 
    \begin{displaymath}
      \Gamma \entailsRanked{\alpha_1}{\zeta}{\rho} A, \Delta\quad
      \text{and} \quad\Gamma', A
      \entailsRanked{\alpha_2}{\zeta}{\rho} \Delta'
    \end{displaymath}
    with  $\alpha_1, \alpha_2 < \alpha$ and $\rank(A) < \rho$,
    then 
    \begin{displaymath}
      \Gamma, \Gamma' \entailsRanked{\alpha}{\zeta}{\rho} \Delta, 
      \Delta'.
    \end{displaymath}
  \item\label{case:comparison1}
    If $\Gamma \entailsRanked{\alpha}{\zeta}{\rho} \Delta$, and if
    $\alpha \leq \beta, \rho\leq \sigma$ with $\beta$ and $\sigma$
    finite, then
    \begin{displaymath}
      \Gamma \entailsRanked{\beta}{\zeta}{\sigma} \Delta
    \end{displaymath}
  \item\label{case:comparison2}
    If $\Gamma \entailsRanked{\alpha}{\zeta}{\rho} \Delta$,
    and if $\zeta < \eta$, then, for any finite $\beta$ and $\sigma$,
    \begin{displaymath}
      \Gamma \entailsRanked{\beta}{\eta}{\sigma} \Delta
    \end{displaymath}.
  \end{enumerate}
\end{definition}
The following lemmas follow by a trivial induction:
\begin{lemma}
  If $\Gamma \entailsMod \Delta$, then, for some finite $\alpha$ and
  $\rho$, and some ordinal $\zeta$,
  \begin{displaymath}
    \Gamma \entailsRanked{\alpha}{\zeta}{\rho}\Delta. 
  \end{displaymath}
\end{lemma}
\begin{lemma}
  A proof $\Gamma \entailsMod \Delta$ is cut free iff it has an
  annotation of the form
  \begin{displaymath}
    \Gamma \entailsRanked{\alpha}{0}{0} \Delta.
  \end{displaymath}
\end{lemma}
\begin{theorem}\label{thm:cutElim}
  Given a McCain-Turner causal theory, the corresponding modal system
  satisfies cut elimination.
\end{theorem}
We prove this by the following lemmas: they are very much the same as
the corresponding lemmas in
\citend[pp.~60ff.]{pohlers89:_proof_theor}. Note that $\alpha\#\beta$
stands for the so-called \emph{natural sum} of the ordinals $\alpha$
and $\beta$ (\citend[p.~43]{pohlers89:_proof_theor}): the relevant
facts that we shall use are that it is defined for any pair of ordinals
$\alpha$ and $\beta$, and that it is strictly monotonic in \emph{both}
arguments.

\begin{lemma}\label{lemma:basicInduction}
  If 
  \begin{displaymath}
    \Gamma \entailsRanked{\alpha}{\zeta}{\rho} X, \Delta
    \quad\text{and}\quad
    \Gamma', X \entailsRanked{\beta}{\eta}{\rho} \Delta'
  \end{displaymath}
  and if $\rank(X) = \rho$, then
  \begin{displaymath}
    \Gamma, \Gamma' \entailsRanked{\alpha + \beta}
    {\zeta\#\eta}{\rho}
    \Delta, \Delta'
  \end{displaymath}
\end{lemma}
\begin{proof}
  We prove this by induction, first on $\zeta \# \eta$, then on
  $\alpha+\beta$. There are four cases.
  \begin{enumerate}
  \item One of the premises is an axiom, $\lfalse\text{L}$, or
    $\ltrue\text{R}$, and $X$ is not principal in it: in this case,
    $\Gamma, \Gamma' \entailsMod \Delta, \Delta'$ is an instance of the
    the same rule, and, by the definition of $\Gamma, \Gamma'
    \entailsRanked{\alpha+\beta}{\zeta\#\eta}{\rho} 
    \Delta, \Delta'$, we have the first case of the result. 
  \item $X$ is non-principal in one or the other of the premises
    (suppose wlog the left one), and that premise is the conclusion of
    a finitary rule. Then the prooftree on the left looks like
    \begin{displaymath}
      \begin{prooftree}
        \begin{Bmatrix}
          \begin{prooftree}
            \Pi_i
            \leadsto
            \Gamma_i \entailsRanked{\alpha_i}{\zeta}{\rho} X, 
            \Delta_i
          \end{prooftree}
        \end{Bmatrix}_{i \in I}
        \justifies
        \Gamma \entailsRanked{\alpha}{\zeta}{\rho} X, \Delta
      \end{prooftree}
    \end{displaymath}
    Since, for all $i$, $\alpha_i < \alpha$ 
    (and, consequently, $\alpha_i + \beta < \alpha + \beta$),    
    we can assume the hypothesis inductively for cuts
    between $\Gamma_i \entailsRanked{\alpha_i}{\zeta_1}{\rho} X, 
    \Delta_i$ and $\Gamma', X \entailsRanked{\beta}{\zeta_2}{\rho}
    \Delta'$, obtaining proofs $\Pi_i$: we complete the proof as follows:
    \begin{displaymath}
      \begin{prooftree}
        \begin{Bmatrix}
          \begin{prooftree}
            \tilde\Pi_i
            \leadsto
            \Gamma_i, \Gamma'\entailsRanked{\alpha_i+\beta}{\zeta\#\eta}
            {\rho}  
            \Delta_i, \Delta'
          \end{prooftree}
        \end{Bmatrix}_{i \in I}
        \justifies
        \Gamma, \Gamma' \entailsRanked{\alpha+\beta}{\zeta\#\eta}{\rho}
        \Delta, \Delta'
      \end{prooftree}
    \end{displaymath}
  \item $X$ is non-principal (wlog on the left), the bottom sequent on
    the left is the conclusion of $\nec\text{R}$, all of the
    premises on the left have $\zeta = \rho = 0$, and the
    sequent on the right has $\eta=\rho=0$. So, the prooftree on
    the left looks like
    \begin{displaymath}
      \begin{prooftree}
        \begin{Bmatrix}
          \begin{prooftree}
            \Pi_i
            \leadsto
            \Gamma_i \entailsRanked{\alpha_i}{0}{0} X, 
            \Delta_i
          \end{prooftree}
        \end{Bmatrix}_{i \in I}
        \justifies
        \Gamma \entailsRanked{\alpha}{0}{0} X, \Delta
      \end{prooftree}.
    \end{displaymath}
    We can apply the
    lemma inductively, obtaining proofs of 
    \begin{displaymath}
      \Gamma_i, \Gamma' \entailsRanked{\alpha_i + \beta}{0}{0}
    \Delta_i, \Delta';
    \end{displaymath}
  we now complete the proof as
    follows
    \begin{displaymath}
      \begin{prooftree}
        \begin{Bmatrix}
          \begin{prooftree}
            \tilde\Pi_i
            \leadsto
            \Gamma_i,\Gamma' \entailsRanked{\alpha_i + \beta}{0}{0}  
            \Delta_i,\Delta'
          \end{prooftree}
        \end{Bmatrix}_{i \in I}
        \justifies
        \Gamma, \Gamma' \entailsRanked{0}{0}{0}  \Delta, \Delta'
      \end{prooftree}
    \end{displaymath}
 \item $X$ is non-principal (wlog on the left), 
    the bottom sequent on the left is the conclusion of
    $\nec\text{R}$, and we do not have
    the previous case. So, the prooftree on the left looks like
    \begin{displaymath}
      \begin{prooftree}
        \begin{Bmatrix}
          \begin{prooftree}
            \Pi_i
            \leadsto
            \Gamma_i \entailsRanked{\alpha_i}{\zeta_i}{\rho_i} X, 
            \Delta_i
          \end{prooftree}
        \end{Bmatrix}_{i \in I}
        \justifies
        \Gamma \entailsRanked{\alpha}{\zeta}{\rho} X, \Delta
      \end{prooftree}
    \end{displaymath}
    where, for all $i$, $\zeta_i < \zeta$. We can apply the
    lemma inductively, obtaining proofs of 
    $\Gamma_i, \Gamma' \entailsRanked{\alpha_i+ \beta}{\zeta_i+ \eta}{\rho'_i}
    \Delta_i, \Delta'$, where, for each $i$,
    $\zeta'_i \leq \zeta_i < \zeta$, and where
    $\rho'_i = \max(\rho_i,\eta)$: we now complete the proof as
    follows
    \begin{displaymath}
      \begin{prooftree}
        \begin{Bmatrix}
          \begin{prooftree}
            \tilde\Pi_i
            \leadsto
            \Gamma_i,\Gamma' \entailsRanked{\alpha'_i}{\zeta'_i}{\rho'_i}  
            \Delta_i,\Delta'
          \end{prooftree}
        \end{Bmatrix}_{i \in I}
        \justifies
        \Gamma, \Gamma' \entailsRanked{0}{\zeta}{0}  \Delta, \Delta'
      \end{prooftree}
    \end{displaymath}
    and apply case~\ref{case:comparison1} of Definition~\ref{definition:entailsRanked}
    to adjust the values of $\alpha$ and $\rho$. 
    
  \item $X$ is principal on both sides; there
    are various cases, depending on
    the principal connective. We treat two 
    of these cases: the others are very similar.
    \begin{description}
    \item[$\rightarrow$]
      The proofs are
      \begin{displaymath}
        \begin{prooftree}
          \begin{prooftree}
            \Pi_1
            \leadsto
            \Gamma, Y \entailsRanked{\alpha_1}{\zeta}{\rho} Z, \Delta
          \end{prooftree}
          \justifies
          \Gamma \entailsRanked{\alpha}{\zeta}{\rho} Y \rightarrow Z, \Delta
          \using
          \rightarrow\text{R}
        \end{prooftree}
        \quad\text{and}\quad
        \begin{prooftree}
          \begin{prooftree}
            \Pi_2
            \leadsto
            \Gamma' \entailsRanked{\beta_1}{\eta}{\rho} Y, \Delta'
          \end{prooftree}
          \quad
          \begin{prooftree}
            \Pi_3
            \leadsto
            \Gamma', Z \entailsRanked{\beta_2}{\eta}{\rho} \Delta'
          \end{prooftree}
          \justifies
          \Gamma', Y \rightarrow Z \entailsRanked{\beta}{\eta}{\rho} \Delta'
          \using
          \rightarrow\text{L}
        \end{prooftree}
      \end{displaymath}
      with $\alpha_1 < \alpha$ and $\beta_1, \beta_2 < \beta$. We
      transform this into
      \begin{displaymath}
        \begin{prooftree}
          \begin{prooftree}
            \tilde\Pi
            \leadsto
            \Gamma, \Gamma' \entailsRanked{\alpha_1 + \beta_1}
            {\zeta\#\eta}{\rho}
            Z, \Delta, \Delta'
          \end{prooftree}
          \quad
          \begin{prooftree}
            \Pi_3
            \leadsto 
            \Gamma', Z \entailsRanked{\beta_2}{\eta}{\rho} \Delta'
          \end{prooftree}
          \justifies
          \Gamma,\Gamma' \entailsRanked{\alpha+\beta}{\zeta\#\eta}{\rho}
          \Delta,\Delta'
        \end{prooftree}
      \end{displaymath}
      where $\tilde\Pi$ is produced from $\Pi_1$ and $\Pi_2$
      by the inductive hypothesis: we have $\rank (Y) , \rank(Z) < \rho$,
      $\alpha_1\#\beta_1 < \alpha\#\beta$, and $\beta_1 < \alpha\#\beta$,
      so the final inference is justified. 
    \item[$\nec$] The proofs are
      \begin{displaymath}
        \begin{prooftree}
          \begin{prooftree}
            \Pi
            \leadsto
            \Gamma \entailsRanked{\alpha_1}{\zeta}{\rho} 
            \phi_1 \land \ldots \land \phi_k, \Delta
          \end{prooftree}
          \justifies
          \Gamma \entailsRanked{\alpha}{\zeta}{\rho} \nec Y, \Delta
          \using
          \nec\text{R}
        \end{prooftree}
        \quad \text{and}\quad
        \begin{prooftree}
          \begin{Bmatrix}
            \begin{prooftree}
              \Pi'_i
              \leadsto
              \Gamma', \phi_{i_1}, \ldots, \phi_{i_{k_i}}
              \entailsRanked{\beta_i}{\eta_i}{\rho} \Delta'
            \end{prooftree}
          \end{Bmatrix}_{i \in I}
          \justifies
          \Gamma', \nec Y \entailsRanked{\beta}{\eta}{\rho} \Delta'
          \using
          \nec\text{L}
        \end{prooftree}
      \end{displaymath}
      where $\alpha_1 < \alpha$, $\eta_i  < \eta$ for all
      $i$.  

      We construct a new proof as follows: the tuple 
      $\phi_1, \ldots, \phi_k$ must match one of the
      tuples indexed by $I$: say it matches $i_0$. 
      So we now have proofs
      \begin{displaymath}
        \begin{prooftree}
          \Pi
          \leadsto
          \Gamma \entailsRanked{\alpha_1}{\zeta}{\rho} 
          \phi_1 \land \ldots \land \phi_k, \Delta
        \end{prooftree}
        \quad \text{and}\quad
        \begin{prooftree}
          \Pi'_{i_0}
          \leadsto
          \begin{prooftree}
            \Gamma', \phi_1, \ldots, \phi_k
            \entailsRanked{\beta_{i_0}}{\eta_{i_0}}{\rho} \Delta'
            \justifies
            \begin{prooftree}
              \Gamma', \phi_1\land \phi_2, \ldots, \phi_k
              \entailsRanked{\beta_{i_0} + 1}{\eta_{i_0}}{\rho} \Delta'
              \leadsto
              \Gamma', \phi_1 \land \ldots l\land \phi_k
              \entailsRanked{\beta_{i_0} + k - 1}{\eta_{i_0}}{\rho}
              \Delta'
              \using
              \land\text{L}
            \end{prooftree}
            \using
            \land\text{L}
          \end{prooftree}
        \end{prooftree}
      \end{displaymath}
      We now have a proof of
      \begin{displaymath}
        \Gamma, \Gamma' \entailsRanked{\alpha'}{\eta_{i_0}}{\rho'} 
        \Delta, \Delta',
      \end{displaymath}
      and $\eta_{i_0} < \eta$, so we can now apply case~\ref{case:comparison2}
      of Definition~\ref{definition:entailsRanked}, and obtain the result.
    \end{description}
  \end{enumerate}
\end{proof}
\begin{proof}[Proof of Theorem~\ref{thm:cutElim}]
  Suppose that we have a proof involving cuts, i.e. a proof of 
  a sequent 
  \begin{equation}
    \Gamma \entailsRanked{\alpha}{\zeta}{\rho+1} \Delta\label{eq:sequent}
  \end{equation}
  We prove that this proof can be replaced by a cut free proof of the
  same sequent by an induction. The inductive hypothesis will be
  that, given such a proof, there is a cut free proof 
  \begin{displaymath}
    \Gamma \entailsRanked{\alpha'}{0}{0} \Delta
  \end{displaymath}
  with some $\alpha'$. We perform the
  induction first on $\zeta$, then on $\rho$, then on $\alpha$. There
  are four cases.
  \begin{enumerate}
  \item If the last inference is an application of $\nec\text{L}$,
    then the premises are of the form 
    $\Gamma_i \entailsRanked{\alpha'_i}{\zeta_i}{\rho_i}\Delta_i$
    with, for all $i$, $\zeta_i < \zeta$. We can assume the
    result inductively, obtaining cut free proofs
    of the premises, and we can complete the
    proof as follows:
    \begin{displaymath}
      \begin{prooftree}
        \{\Gamma_i \entailsRanked{\alpha''_i}{0}{0}
          \Delta_i\}_{i \in I}
        \justifies
        \Gamma \entailsRanked{0}{0}{0} \Delta
        \using \nec\text{L}
      \end{prooftree},
    \end{displaymath}
    which is a cut free proof of \eqref{eq:sequent}. 
  \item If the last inference is a cut of rank less than $\rho$,
    then the premises will have the same value of $\zeta$
    and $\rho$, but  smaller values of $\alpha$: inductively,
    we can assume that these premises have cut free proofs. 
    So we now have a proof as follows:
    \begin{displaymath}
      \begin{prooftree}
        \Gamma' \entailsRanked{\alpha'}{0}{0} X, \Delta'
        \quad
        \Gamma'',X \entailsRanked{\alpha''}{0}{0}\Delta''
        \justifies
        \Gamma', \Gamma'' \entailsRanked{0}{0}{\rho'} \Delta', \Delta''
      \end{prooftree}
     \end{displaymath}
     with $\rho' < \rho + 1$. So, by our inductive hypothesis, we can
     assume the result.
  \item If the last inference is a cut of rank $\rho$, then its
    premises will be of the form 
    \begin{align*}
      \Gamma'
      &\entailsRanked{\alpha'}{\zeta}{\rho+1}\Delta'\quad\text{and}\quad
      \Gamma''\entailsRanked{\alpha''}{\zeta}{\rho+1}\Delta''
      \intertext{We can apply the result inductively to obtain cut
        free proofs of the premises, which (after applying
        Definition~\ref{definition:entailsRanked} 
        case~\ref{case:comparison1}) can be assumed to be of the form}
      \Gamma' &\entailsRanked{\alpha'}{0}{\rho}\Delta'\quad\text{and}\quad
      \Gamma''\entailsRanked{\alpha''}{0}{\rho}\Delta'';
    \end{align*}  
    we then apply Lemma~\ref{lemma:basicInduction} to obtain a proof of 
    \begin{displaymath}
      \Gamma \entailsRanked{\alpha + \alpha'}{0}{\rho} \Delta. 
    \end{displaymath}
    We have possibly reduced the value of $\zeta$, and
    certainly reduced the value of $\rho$, so we can assume
    the result by the inductive hypothesis. 
  \item If the last inference is an application of a non-cut finitary
    rule, then its premises must have the same value of $\zeta$ and
    $\rho$, but smaller values of $\alpha$, than the conclusion:
    we can again apply the result inductively to the premises,
    and then apply the rule to the premises, giving a cut free
    proof of \eqref{eq:sequent}.
  \end{enumerate}
\end{proof}
\section{Turner's 
  \emph{Logic of Universal Causation}}\label{sec:turn-emphl-univ}
As we have indicated, Turner~\cite{turner99:_theor_univer_causat}
also has a modal system which he uses as a metatheory for the
McCain-Turner procedure. His treatment has the
following features:
\begin{enumerate}
\item We are given an \textbf{S5} modal operator, written \turnerMod,
  and a theory $T$ in the language of that operator
\item Recall that a Kripke model of \textbf{S5} is
  \begin{enumerate}
  \item a set $\mathcal{W}$ of worlds,
  \item a truth-functional forcing relation between
    worlds and propositions of the non-modal language, and
  \item an equivalence relation on the set of worlds. 
  \end{enumerate}

\item[]  We then say that an \textbf{S5} Kripke model of the theory
  $T$ is \emph{causally explained}
  if 
  \begin{enumerate}
    \item it has a single world $w$
    \item if any model of $T$ has a set of worlds $\mathcal{W}$
      which is a superset of $\{w\}$, with a forcing relation
      extending that for $w$, then $\{w\}$ is a single
      equivalence class in the larger model. 
      
      (Equivalently, the inclusion of $w$ in $\mathcal{W}$ is what is
      called a \emph{p-morphism}; see
      \citend{vanBenthemJ:hanplii,benthem98:_manual_inten_logic}).
  \end{enumerate}
\item Then we say that a proposition is \emph{causally explained} if 
  it is forced at the unique world in all causally explained models. 
\end{enumerate}
Now this definition of a causally explained model involves not simply
a model, but all models which contain it (in some appropriate sense of
`contain'): it is not obvious that this definition can be replaced
with one which only talks, in the standard way, about theoremhood in a
single model. In fact, as we shall see, it cannot.

\begin{example}\label{xmpl:TurnerNotMono}
  Consider a language with a single atom $p$,
  and two theories in that language:
  \begin{align}
    T_1 & = \{p \rightarrow \turnerMod p \}\\
    T_2 & = \{p \rightarrow \turnerMod p,\lnot p \rightarrow
    \turnerMod \lnot p\}
  \end{align}
  Then $T_1$ has a single causally explained model: its world
  forces $p$. $T_2$ has two causally explained models, one with
  a single world forcing $p$, and one with a single world
  forcing $\lnot p$. Consequently, the proposition $p$ is
  causally explained by $T_1$, but not by $T_2$.
\end{example}
This has, as a consequence,
\begin{proposition}\label{proposition:TurnerNotMono}
  There is no set of propositions $\Gamma$ such that,
  for any \textbf{S5} theory $T$,
  \begin{displaymath}
    p \text{ is causally explained by } T 
    \quad\Leftrightarrow\quad
    \Gamma, T \entails_{\textbf{S5}} p.
  \end{displaymath}
\end{proposition}
\begin{proof}
  If there were such a set of propositions $\Gamma$, then
  the relation
  \begin{displaymath}
    p \text{ is causally explained by }T
  \end{displaymath}
  would be monotonic in $T$: as Example~\ref{xmpl:TurnerNotMono}
  shows, it is not. 
\end{proof}

By contrast, our logic is monotonic: it is given by a sequent calculus
of the normal sort, \emph{with} the weakening rule. 

\begin{example}
  We should consider the analogue of Example~\ref{xmpl:TurnerNotMono}
  in our system. The two theories would correspond to two different
  modal operators, $\nec_1$ and $\nec_2$, and we have :
  \begin{align*}
    \nec_1 p \lequiv p &\quad\nec_1 \lnot p \lequiv \lfalse\\
    \nec_2 p \lequiv p &\quad\nec_2 \lnot p \lequiv \lnot p
  \end{align*}
  The two sets of propositions $\Gamma_i$ are as follows:
  \begin{align*}
    \Gamma_1 & = \{ \nec_1 p \leftrightarrow p, 
    \nec_1\lnot p \leftrightarrow \lnot p\}\\
    & = \{ p \leftrightarrow p, \lnot \lfalse \leftrightarrow \lnot p \}\\
    & = \{ \ltrue, p \}\\
    \Gamma_2 & = \{ \nec_2 p \leftrightarrow p, 
    \nec_2\lnot p \leftrightarrow \lnot p\}\\
    & = \{ p \leftrightarrow p, \lnot p \leftrightarrow \lnot p \}\\
    & = \{ \ltrue\}.
  \end{align*}
  This is in accordance with the expected results: for $T_1$, only
  the world which forces $p$ is causally explained, and so $\Gamma_1$
  (essentially) contains only $p$. For $T_2$, both
  worlds are causally explained (both the one which forces $p$
  and the one which forces $\lnot p$), and $\Gamma_2$ is essentially
  vacuous. 
\end{example}
So, although our modal logic and Turner's yield the same sets of
causally explained models (and hence the same sets of causally
explained propositions), they have radically different properties.
Ours is monotonic and has a well-behaved proof theory: Turner's 
is non-monotonic and, as it stands, does not have a proof theory.
Although Turner's system is formulated using
\textbf{S5}, Proposition~\ref{proposition:TurnerNotMono} shows that
we cannot, in any obvious way, use the proof theory of \textbf{S5}
as a proof theory for his system. 

\bibliographystyle{jflnat}
\bibliography{frame,WhiteGG,semantics,ai,phil,prooftheory,dynamics,mas,logic,ceb,semantics}
\end{document}